\documentclass[a4paper,11pt]{amsart}
\usepackage{mathrsfs}
\usepackage{cases}
\usepackage{amsfonts}
\usepackage{graphicx}
\usepackage{amsmath}
\usepackage{amssymb}
\usepackage[all]{xypic}
\usepackage[all]{xy}
\usepackage{color}
\usepackage{colordvi}
\usepackage{multicol}
\usepackage{polynom}
\usepackage[colorlinks,citecolor=black,linkcolor=black]{hyperref}
\begin{document}
\input xy
\xyoption{all}

\newcommand{\ind}{\operatorname{inj.dim}\nolimits}
\newcommand{\id}{\operatorname{id}\nolimits}
\newcommand{\Mod}{\operatorname{Mod}\nolimits}
\newcommand{\End}{\operatorname{End}\nolimits}
\newcommand{\Ext}{\operatorname{Ext}\nolimits}
\newcommand{\Hom}{\operatorname{Hom}\nolimits}
\newcommand{\aut}{\operatorname{Aut}\nolimits}
\newcommand{\Ker}{{\operatorname{Ker}\nolimits}}
\newcommand{\Iso}{\operatorname{Iso}\nolimits}
\newcommand{\Coker}{\operatorname{Coker}\nolimits}
\renewcommand{\dim}{\operatorname{dim}\nolimits}
\newcommand{\Cone}{{\operatorname{Cone}\nolimits}}
\renewcommand{\Im}{\operatorname{Im}\nolimits}

\newcommand{\cc}{{\mathcal C}}
\newcommand{\ce}{{\mathcal E}}
\newcommand{\cs}{{\mathcal S}}
\newcommand{\cf}{{\mathcal F}}
\newcommand{\cx}{{\mathcal X}}
\newcommand{\cy}{{\mathcal Y}}
\newcommand{\cl}{{\mathcal L}}
\newcommand{\ct}{{\mathcal T}}
\newcommand{\cu}{{\mathcal U}}
\newcommand{\cm}{{\mathcal M}}
\newcommand{\cv}{{\mathcal V}}
\newcommand{\ch}{{\mathcal H}}
\newcommand{\ca}{{\mathcal A}}
\newcommand{\mcr}{{\mathcal R}}
\newcommand{\cb}{{\mathcal B}}
\newcommand{\ci}{{\mathcal I}}
\newcommand{\cj}{{\mathcal J}}
\newcommand{\cp}{{\mathcal P}}
\newcommand{\cg}{{\mathcal G}}
\newcommand{\cw}{{\mathcal W}}
\newcommand{\co}{{\mathcal O}}
\newcommand{\cd}{{\mathcal D}}
\newcommand{\cn}{{\mathcal N}}
\newcommand{\ck}{{\mathcal K}}
\newcommand{\calr}{{\mathcal R}}
\newcommand{\ol}{\overline}
\newcommand{\ul}{\underline}
\newcommand{\cz}{{\mathcal Z}}
\newcommand{\st}{[1]}
\newcommand{\ow}{\widetilde}
\renewcommand{\P}{\mathbf{P}}
\newcommand{\pic}{\operatorname{Pic}\nolimits}
\newcommand{\Spec}{\operatorname{Spec}\nolimits}
\newtheorem{theorem}{Theorem}[section]
\newtheorem{acknowledgement}[theorem]{Acknowledgement}
\newtheorem{algorithm}[theorem]{Algorithm}
\newtheorem{axiom}[theorem]{Axiom}
\newtheorem{case}[theorem]{Case}
\newtheorem{claim}[theorem]{Claim}
\newtheorem{conclusion}[theorem]{Conclusion}
\newtheorem{condition}[theorem]{Condition}
\newtheorem{conjecture}[theorem]{Conjecture}
\newtheorem{construction}[theorem]{Construction}
\newtheorem{corollary}[theorem]{Corollary}
\newtheorem{criterion}[theorem]{Criterion}
\newtheorem{definition}[theorem]{Definition}
\newtheorem{example}[theorem]{Example}
\newtheorem{exercise}[theorem]{Exercise}
\newtheorem{lemma}[theorem]{Lemma}
\newtheorem{notation}[theorem]{Notation}
\newtheorem{problem}[theorem]{Problem}
\newtheorem{proposition}[theorem]{Proposition}
\newtheorem{remark}[theorem]{Remark}
\newtheorem{solution}[theorem]{Solution}
\newtheorem{summary}[theorem]{Summary}
\newtheorem*{thm}{Theorem}

\renewcommand{\theequation}{\arabic{section}.\arabic{equation}}

\def \bp{{\mathbf p}}
\def \bA{{\mathbf A}}
\def \bL{{\mathbf L}}
\def \bF{{\mathbf F}}
\def \bS{{\mathbf S}}
\def \bC{{\mathbf C}}
\def \bD{{\mathbf D}}
\def \Z{{\Bbb Z}}
\def \F{{\Bbb F}}
\def \C{{\Bbb C}}
\def \N{{\Bbb N}}
\def \Q{{\Bbb Q}}
\def \G{{\Bbb G}}
\def \X{{\Bbb X}}
\def \P{{\Bbb P}}
\def \K{{\Bbb K}}
\def \E{{\Bbb E}}
\def \A{{\Bbb A}}
\def \BH{{\Bbb H}}
\def \T{{\Bbb T}}
\def\rH{{\mathrm{H}}}

\title[Green's formula and Derived Hall algebras]{From Green's formula to Derived Hall algebras}
{\thanks{The author is supported by the National Natural Science Foundation of China (No. 12001107), the Discipline (Profession) Leader Training Project of Anhui Province(No. DTR2024026)
, the Scientific Research and Innovation Team Project of Fuyang Normal University (No. TDJC2021009) and the University Outstanding Youth Research Project of Anhui Province(No. 2022AH020082)}}
\author[Ji Lin]{Ji Lin}
\address{Department of Mathematics and Statistics, Fuyang Normal University, Fuyang 236037, P.R.China}
\email{jlin@fynu.edu.cn}


\begin{abstract}
The aim of this note is to clarify the relationship between Green's formula and the associativity of multiplication for derived Hall algebra in the sense of To\"{e}n (Duke Math J 135(3):587-615, 2006), Xiao and Xu (Duke Math J 143(2):357-373, 2008) and Xu and Chen (Algebr Represent Theory 16(3):673-687, 2013). Let $\mathcal{A}$ be a finitary hereditary abelian category. It is known that the associativity of derived Hall algebra $\mathcal{D}\mathcal{H}_t(\mathcal{A})$ implies Green's formula. We show the converse statement holds. Namely, Green's formula implies the associativity of the derived Hall algebra $\mathcal{D}\mathcal{H}_t(\mathcal{A})$. \end{abstract}

\maketitle
%\tableofcontents
\section{Introduction}
Ringel-Hall algebra is an important tool to encode the extension structures of categories. By definition, the Ringel-Hall algebra $\mathcal{H}(\mathcal{A})$ associated to a finitary abelian category $\mathcal{A}$ is a vector space over $\mathbb{Q}$ (or $\mathbb{C}$) generated by the isomorphism classes $[X]$ of $\mathcal{A}$ with the multiplication
$$[X]\diamond[Y]=\sum_{[L]\in\Iso(\mathcal{A})}g_{XY}^L[L]$$
for any $[X],[Y]\in\Iso(\mathcal{A})$, where $g_{XY}^L$ called the {\it Hall number} is the number of subobjects $L'$ of $L$ such that $L'\cong Y$ and $L/L'\cong X$. For a hereditary algebra $A$ over a finite field $\mathbb{F}_q$, C. M. Ringel and J. A. Green showed that the positive part of the quantum group of the corresponding Kac-Moody algebra can be realized as the composition subalgebra of the Ringel-Hall algebra $\mathcal{H}(A)$ associated to the module category of $A$. In particular, Green also introduced a bialgebra structure of $\mathcal{H}(A)$, of which the so-called Green's formula is the key ingredient. Ringel asked that whether one can recover the whole quantum group by Ringel-Hall algebras. A direct way is to glue together two Borel parts by using reduced Drinfeld double as shown in \cite{X97}. Ringel-Hall algebras have been intensively studied in the past three decades, and we refer to \cite{R90, Rie94, Gr95,R96, DX04, Br13} for more details.

To realize the quantum group intrinsically one need to replace the category by a larger one. In \cite{T06} To\"{e}n made a remarkable development in this direction. He constructed a derived Hall algebra over a dg-category satisfying some finiteness conditions which is closely related to Ringel-Hall algebra in the case of hereditary abelian categories. By counting triangles and applying the octahedral axiom, Xiao and Xu \cite{XX08} found a direct way to prove the formula of To\"{e}n for a triangulated categories satisfying the left homological finiteness conditions. Later Xu and Chen \cite{XC13} showed an analogue of To\"{e}n's formula in the case of triangulated categories with odd periodic translation functors.  In particular, the bounded derived category $\mathcal{D}^b(\mathcal{A})$ (setting $\mathcal{D}_0(\mathcal{A})=\mathcal{D}^b(\mathcal{A})$ by convention) of $\mathcal{A}$ satisfies the left homological finiteness conditions required in \cite{T06, XX08} and the relative derived category of $\mathcal{A}$ denoted by $\mathcal{D}_t(\mathcal{A})$ satisfies the homological finiteness conditions in \cite{XC13} for an odd positive integer $t$. However it fails to define Hall algebras over $2$-periodic triangulated categories and one cannot recover the whole quantum groups via derived Hall algebras.

T. Bridgeland considered the Ringel-Hall algebra of the exact category of the $2$-periodic complexes with projective components over an abelian category which has enough projectives, and he realized the whole quantum group via the so-called Bridgeland's Hall algebra which is, by definition, the localization of Ringel-Hall algebra of the $2$-periodic complexes at the set of acyclic complexes. Later on, M. Gorsky \cite{Gor13, Gor16} constructed the semi-derived Ringel-Hall algebras over the categories of complexes over exact categories and Frobenius categories respectively, which is, in some sense, between the usual Hall algebras and the Hall algebras of the bounded derived categories.  However in the case of $2$-periodic complexes they both requires the condition that the exact categories have enough projective objects. And in this case the both algebras can be used to realize the whole quantum group.

Motivated by the work of Bridgeland and Gorsky, Lu and Peng \cite{LuP21} introduced the semi-derived Ringel-Hall algebra over a hereditary abelian category which may not have enough projective objects. And it is called {\it modified Ringel-Hall algebra} in the previous version of their paper \cite{LuP16} and in our previous work \cite{LinP2019}. We studied the semi-derived Ringel-Hall algebra in the case of bounded complexes in \cite{LinP2019}, we also found that Green's formula can be deduced from the associativity of semi-derived Ringel-Hall algebra. Similarly one can obtain Green's formula by the associativity of derived Hall algebras as shown in \cite{Zhang, Ru21}. On the contrary, we mainly consider the opposite problem that if the derived Hall algebras associated to the categories of the complexes over hereditary abelian categories can be deduced by Green's formula.

Most recently, H. Zhang \cite{Zhang22} defined a Hall algebra for the root category by applying the derived Hall numbers of bounded derived category. By the same way, he also gave an explicit formula for odd-periodic derived Hall algebras in \cite{Zhang23}. The main tools he used are Green's formula and the relation between the five-term exact sequences and the extensions in bounded derived categories. On the other side, J. Chen, M. Lu and S. Ruan \cite{CLR23} defined a Hall algebra for the root category directly by counting the triangles in the root category. They \cite{CLR22} also used the derived Hall algebra of 1-periodic complexes to realize the $\imath$quantum groups of split type. Of course, the structure constant given in this note is the same as that in \cite{Zhang23} for the case of odd-periodic complexes, also equivalent to  that in \cite{CLR22} for the case of 1-periodic complexes. However, we gave direct characterizations of derived Hall numbers for both bounded derived categories and odd-periodic derived categories by using the classic Hall numbers of abelian categories, different to \cite{Zhang22, CLR22,CLR23}. Furthermore, we directly showed the associativity of derived Hall algebras by Green's formula.    

The paper is organized as follows. In Section \ref{Sec Preliminaries}, we collect some definitions and results on relative derived categories and Green's formula. In Section \ref{Sec ass algebras}, by using Green's formula, we directly define an associative algebra $\cl_t(\ca)$ over the derived category $\cd_t(\ca)$. Section \ref{Sec derived Hall algebras} is devoted to recalling the definitions, generators and generating relations of derived Hall algebras in the case of bounded complexes and odd periodic complexes of the hereditary abelian category $\mathcal{A}$ respectively. In Section \ref{Sec derived Hall numbers}, we first work out the structure constants of derived Hall algebras following the methods of Y. Chen and F. Wu in \cite{Chen19, W18} and then show our main result by the natural isomorphism between $\cl_t(\ca)$ and derived Hall algebra $\cd_t(\ca)$ which states that Green's formula implies the associativity of derived Hall algebras in both cases of bounded complexes and odd periodic complexes of a finitary hereditary abelian category $\mathcal{A}$.

Throughout the paper we denote by $|S|$ the cardinality of a finite set $S$. Without specification, $\mathcal{A}$ is a finitary abelian category over a finite field $k=\mathbb{F}_q$, i.e.,
a small abelian cateory satisfying: $(1)~|\Hom_{\mathcal{A}}(X, Y)|<\infty$ and $(2)~|\Ext^1_{\mathcal{A}}(X, Y)|<\infty$ for any $X, Y\in\mathcal{A}$. For any $X$ in a category, we always use $[X]$ to denote the isomorphism class of $X$, and denote by $\aut(X)$ the auto morphism group of $X$ and $a_X$ the cardinality of $\aut(X)$. And $K_0(\mathcal{A})$ is the Grothendieck group of $\mathcal{A}$ with its elements denoted by $\widehat{X}$ for all $X\in\mathcal{A}$.

\vspace{0.2cm} \noindent{\bf Acknowledgments.}
The author is deeply indebted to Professors Changjian Fu and Ming Lu for their helpful discussions and comments.

\section{Preliminaries}\label{Sec Preliminaries}
In this section, we recall relative derived categories of $t$-periodic complexes over $\mathcal{A}$ and Green's formula (see \cite{PX97,Br13,Gr95}).
\subsection{Relative derived categories}Let $\mathcal{A}$ be a hereditary abelian category. For each integer $t\geq 1$, write $\mathbb{Z}_t=\mathbb{Z}/t=\{0,1,\cdots,t-1\}$. By definition, a $t$-periodic complex over $\mathcal{A}$ denoted by $M^{\bullet}=(M^i, d^i_{M^{\bullet}})_{i\in\mathbb{Z}_t}$ is composed of objects $M^i\in\mathcal{A}$ and morphisms $d^i_{M^{\bullet}}: M^i\rightarrow M^{i+1}$ for $i\in\mathbb{Z}_t$ satisfying $d^{i+1}_{M^{\bullet}}d^i_{M^{\bullet}}=0$. For convenience, we also express $M^{\bullet}=(M^i, d^i_{M^{\bullet}})_{i\in\mathbb{Z}_t}$ in the following form:
\[M^0\xrightarrow{d^0_{M^{\bullet}}}M^1\xrightarrow{d^1_{M^{\bullet}}}M^2\rightarrow
\cdots\rightarrow M^{t-1}\xrightarrow{d^{t-1}_{M^{\bullet}}}M^0.\]
Let $L^{\bullet}=(L^i, d^i_{L^{\bullet}})_{i\in\mathbb{Z}_t}$ and $M^{\bullet}=(M^i, d^i_{M^{\bullet}})_{i\in\mathbb{Z}_t}$ be two $t$-periodic complexes, a {\it morphism} $$f^\bullet: L^{\bullet}\rightarrow M^{\bullet}$$ is a family of morphisms $f^i: L^i\rightarrow M^i$ satisfying $f^{i+1}d^i_{L^{\bullet}}=d^i_{M^{\bullet}}f^i$ for all $i\in\mathbb{Z}_t$. The $i$th {\it homology} of $M^{\bullet}$ denoted by $\rH^i(M^\bullet)$ is by definition $\Ker(d^i_{M^{\bullet}})/\Im(d^{i-1}_{M^{\bullet}})$ for any $i\in\mathbb{Z}_t$. The {\it category of $t$-periodic complexes} over $\mathcal{A}$ is denoted by $\mathcal{C}_t(\mathcal{A})$.

Let $f^\bullet,g^\bullet: L^{\bullet}\rightarrow M^{\bullet}$ be two morphisms of $t$-periodic complexes. A {\it relative homotopy} $s^\bullet$ from $f^\bullet$ to $g^\bullet$ is a sequence of morphisms $s^i: L^i\rightarrow M^{i-1}$ of $\mathcal{A}$ such that
\[f^i-g^i=d^{i-1}_{M^{\bullet}}s^i+s^{i+1}d^{i}_{L^{\bullet}}\]
for all $i\in\mathbb{Z}_t$. Morphisms $f^\bullet$ and $g^\bullet$ are said to be {\it relatively homotopic} if there exists a relative homotopy from $f^\bullet$ to $g^\bullet$. Relative homotopy is an equivalence relation compatible with composition of morphisms, so one can form the {\it relative homotopy category} $\mathcal{K}_t(\mathcal{A})$ of $t$-periodic complexes over $\mathcal{A}$ by considering all $t$-periodic complexes as objects and the additive group of relative homotopy classes of morphisms from $L^\bullet$ to $M^\bullet$ in $\mathcal{C}_t(\mathcal{A})$ as the group of morphisms from $L^\bullet$ to $M^\bullet$.

As in categories of bounded complexes, one can define quasi-isomorphisms in $\mathcal{C}_t(\mathcal{A})$ and $\mathcal{K}_t(\mathcal{A})$ and similarly form a new additive category $\mathcal{D}_t(\mathcal{A})$ by localizing $\mathcal{K}_t(\mathcal{A})$ with respect to the set of all quasi-isomorphisms, called the {\it relative derived category} of $t$-periodic complexes over $\mathcal{A}$. Similar to the appendix of \cite{PX97}, $\mathcal{K}_t(\mathcal{A})$ and $\mathcal{D}_t(\mathcal{A})$ are both triangulated categories with the translation functor-the shift functor $[1]$.

Simply we set $\mathbb{Z}_0=\mathbb{Z}$ and write $\mathcal{C}_0(\mathcal{A})=\mathcal{C}^b(\mathcal{A})$ and $\mathcal{D}_0(\mathcal{A})=\mathcal{D}^b(\mathcal{A})$. For any $X\in\mathcal{A}$ and $i\in\mathbb{Z}_t$, we simply denote by $Z_X^{[n]}$ the stalk complex with its nonzero component $X$ sitting in the degree $n$. In particular, if $t=1$, set $Z_X:=Z_{X,0}$. And for each $t\geq 0$, we have the following isomorphism.

\begin{proposition}[\cite{LinP2021}, Proposition 2.2]\label{proposition t-periodic iso to the homology}
Let $\mathcal{A}$ be a hereditary abelian category and $M^{\bullet}=(M^i, d^i_{M^{\bullet}})_{i\in\mathbb{Z}_t}\in\cc_t(\ca)$. Then $M^\bullet$ is isomorphic to $\bigoplus\limits_{i\in\mathbb{Z}_t}Z_{\rH^i(M^\bullet)}^{[i]}$ in $\cd_t(\ca)$.
\end{proposition}

\subsection{Green's formula}For any objects $A, B, C\in\ca$,  denote by $g^C_{AB}$ the Hall number. Then one has the following homological formula (cf. \cite{P97, Rie94}) $$g^C_{AB}=\frac{|\Ext_{\ca}^1(A, B)_C|}{|\Hom_{\ca}(A, B)|}\frac{a_C}{a_Aa_B},$$
where $\Ext_{\ca}^1(A, B)_C$ is the subset of $\Ext_{\ca}^1(A, B)$  consisting of the extensions with the middle term isomorphic to $C$. The following is Green's formula.
\begin{theorem}[\cite{Gr95},Theorem 2]\label{Green's formula}
Let $A, B, A', B'$ be fixed objects of $\ca$. Then there holds
\begin{eqnarray*}
&&a_Aa_Ba_{A'}a_{B'}\sum\limits_{C\in\Iso(\ca)}g_{AB}^{C}g_{A'B'}^{C}\frac{1}{a_C}\\
&=&\sum\limits_{X,Y,X',Y'\in\Iso(\ca)}\frac{|\Ext_{\ca}^1(X,Y')|}{|\Hom_{\ca}(X,Y')|}
g_{XX'}^{A}g_{YY'}^{B}g_{XY}^{A'}g_{X'Y'}^{B'}a_Xa_Ya_{X'}a_{Y'}.
\end{eqnarray*}
\end{theorem}
\section{Associative algebra $\cl_t(\ca)$ over derived category}\label{Sec ass algebras}
Assume that $\mathcal{A}$ is a finitary hereditary abelian category over the finite field $k$. We denote by $\langle-,-\rangle$ the multiplicative Euler form defined by $$\langle A, B\rangle:=\frac{|\Hom_{\ca}(A, B)|}{|\Ext^1_{
\ca}(A, B)|}$$
for objects $A, B\in\ca$, and it can be obviously descent to the Grothendieck group $K_0(\ca)$.

For any $t=0$ or an odd positive integer, let $\mathcal{L}_t(\ca)$ be the $\mathbb{C}$-vector space with the basis $\{[A^\bullet]|A^\bullet=(A^i,0)_{i\in\mathbb{Z}_t}\in\cd_{t}(\ca)\}$, and with the multiplications defined on the basis elements by the following rules. Let $A_1^\bullet=(A_1^i, 0)_{i\in\mathbb{Z}}, A_2^\bullet=(A_2^i, 0)_{i\in\mathbb{Z}}\in\cd_0(\ca)$, define the multiplication by the formula
\begin{eqnarray*}
[A_1^\bullet][A_2^\bullet]&=&\sum\limits_{X^{\bullet}=(X^i,0)_{i\in\mathbb{Z}}\in\Iso(\cd_0(\ca))}\sum_{\begin{array}{c}I^i,M^i,N^i \in\Iso(\ca)\\i\in\mathbb{Z}\end{array}}\\
&&\left(\prod_{i\in\mathbb{Z}}g^{A_1^i}_{I^iM^i}g^{X^i}_{M^iN^i}g^{A_2^i}_{N^iI^{i-1}}
\frac{a_{N^i}a_{M^i}a_{I^i}}{a_{A_1^i}a_{A_2^i}}\frac{\prod_{k\geq 2}\langle\widehat{A_2^{i+k}},\widehat{A_1^{i}}\rangle^{(-1)^k}}{\langle\widehat{N^{i+1}},\widehat{M^i}\rangle}\right)[X^\bullet].
\end{eqnarray*}
If $t$ is an odd positive integer, for any $A_1^\bullet=(A_1^i, 0)_{i\in\mathbb{Z}_t}, A_2^\bullet=(A_2^i, 0)_{i\in\mathbb{Z}_t}, X^\bullet=(X^i, 0)_{i\in\mathbb{Z}_t}\in\cd_t(\mathcal{A})$, let
\begin{eqnarray*}
&&h_{A_1^\bullet A_2^\bullet}^{X^\bullet}\\
&=&\sum_{\begin{array}{c}
M^i,N^i,S^i\in\Iso(\mathcal{A})\\i\in \mathbb{Z}_t\end{array}}\left(\prod_{i\in\mathbb{Z}_t}\frac{g_{N^iS^i}^{A_2^i}g_{M^iN^i}^{X^i}g_{S^{i+1}M^i}^{A_1^i}a_{M^i}a_{N^i}a_{S^i}}
{a_{X^i}\langle A_1^i,S^i\rangle\langle S^{i+1},N^i\rangle}\sqrt{\langle A_1^i,A_2^i\rangle\prod_{k=1}^{t-1}\langle A_1^{i+k},A_2^i\rangle^{(-1)^{k+1}}}\right),
\end{eqnarray*}
and the multiplication is defined by
\begin{eqnarray*}
[A_1^\bullet][A_2^\bullet]&=&\sum\limits_{X^{\bullet}=(X^i,0)_{i\in\mathbb{Z}_t}\in\Iso(\cd_t(\ca))}h_{A_1^\bullet A_2^\bullet}^{X^\bullet}\frac{a_{X^\bullet}'}{a_{A_1^\bullet}'
a_{A_2^\bullet}'}[X^\bullet],
\end{eqnarray*}
where $a_{Y^\bullet}'=a_{Y^\bullet}\sqrt{\prod_{i=1}^{t}|\Hom_{\cd_t(\ca)}(Y^\bullet[i], Y^\bullet)|^{(-1)^i}}$ for any $Y^\bullet\in\cd_t(\ca)$.

\begin{theorem}\label{Thm associativity}
For any $t=0$ or an odd positive integer, $\mathcal{L}_t(\ca)$ defined as above are associative algebras with unity $[0]\in\Iso(\cd_t(\ca))$.
\end{theorem}
\begin{proof}
We only need to show the associativity of the multiplication and split the proof into two cases. For $t=0$, let $A_1^\bullet=(A_1^i, 0)_{i\in\mathbb{Z}}, A_2^\bullet=(A_2^i, 0)_{i\in\mathbb{Z}}, A_3^\bullet=(A_3^i, 0)_{i\in\mathbb{Z}}\in\cd_0(\ca)$, by definition
\begin{eqnarray*}
&&\left([A_1^\bullet][A_2^\bullet]\right)[A_3^\bullet]\\
&=&\sum_{\begin{array}{c}X^\bullet\in\Iso(\cd_0(\ca))\end{array}}\\
&&\left(\sum_{\begin{array}{c}M_1^i,N_1^i,I_1^i\\i\in\mathbb{Z}\end{array}}\left(\prod_{i\in\mathbb{Z}}
g^{A_1^i}_{I_1^iM_1^i}g^{X^i}_{M_1^iN_1^i}g^{A_2^i}_{N_1^iI_1^{i-1}}
\frac{a_{N_1^i}a_{M_1^i}a_{I_1^i}}{a_{A_1^i}a_{A_2^i}}\frac{\prod_{k\geq 2}\langle\widehat{A_2^{i+k}},\widehat{A_1^{i}}\rangle^{(-1)^k}}
{\langle\widehat{N_1^{i+1}},\widehat{M_1^i}\rangle}\right)\right)[X^\bullet][A_{3}^\bullet]
\end{eqnarray*}
\begin{eqnarray*}
&=&\sum_{\begin{array}{c}Y^\bullet, X^\bullet\in\Iso(\cd_0(\ca))\end{array}}\sum_{\begin{array}{c}M_1^i,N_1^i,I_1^i\\
M_2^i,N_2^i,I_2^i\\i\in\mathbb{Z}\end{array}}\prod_{i\in\mathbb{Z}}
g^{A_1^i}_{I_1^iM_1^i}g^{X^i}_{M_1^iN_1^i}g^{A_2^i}_{N_1^iI_1^{i-1}}g^{X^i}_{I_2^iM_2^i}
g^{Y^i}_{M_2^iN_2^i}g^{A_3^i}_{N_2^i I_2^{i-1}}
\\
&&\frac{a_{N_1^i}a_{M_1^i}a_{I_1^i}a_{N_2^i}a_{M_2^i}a_{I_2^i}}{a_{A_1^i}a_{A_2^i}a_{A_3^i}a_{X^i}}\frac{\prod_{k\geq 2}\left(\langle\widehat{A_2^{i+k}},\widehat{A_1^{i}}\rangle\langle\widehat{A_3^{i+k}},\widehat{X^{i}}\rangle\right)^{(-1)^k}}
{\langle\widehat{N_1^{i+1}},\widehat{M_1^i}\rangle\langle\widehat{N_2^{i+1}},\widehat{M_2^i}\rangle}[Y^\bullet]\\
&=&\sum_{Y^\bullet\in\Iso(\cd_0(\ca))}\sum_{\begin{array}{c}M_1^i,N_1^i,I_1^i\\M_2^i,N_2^i,I_2^i\\i\in\mathbb{Z}\end{array}}\left(\prod_{i\in\mathbb{Z}}
\left(\sum_{\begin{array}{c}X^i\end{array}}g^{X^i}_{M_1^iN_1^i}g^{X^i}_{I_2^iM_2^i}
\frac{a_{N_1^i}a_{M_1^i}a_{M_2^i}a_{I_2^i}}{a_{X^i}}\right)
\right.\\
&&\left.g^{A_1^i}_{I_1^iM_1^i}g^{A_2^i}_{N_1^iI_1^{i-1}}
g^{A_3^i}_{N_2^i I_2^{i-1}}g^{Y^i}_{M_2^iN_2^i}\frac{a_{I_1^i}a_{N_2^i}}
{a_{A_1^i}a_{A_2^i}a_{A_3^i}}\frac{\prod_{k\geq 2}\left(\langle\widehat{A_2^{i+k}},\widehat{A_1^{i}}\rangle\langle\widehat{A_3^{i+k}},\widehat{X^{i}}\rangle\right)^{(-1)^k}}
{\langle\widehat{N_1^{i+1}},\widehat{M_1^i}\rangle\langle\widehat{N_2^{i+1}},\widehat{M_2^i}\rangle}\right)[Y^\bullet].
\end{eqnarray*}

By Green's formula (see Theorem \ref{Green's formula}), we obtain that
\begin{eqnarray*}
&&\left([A_1^\bullet][A_2^\bullet]\right)[A_3^\bullet]\\
&=&\sum_{Y^\bullet\in\Iso(\cd_0(\ca))}\sum_{\begin{array}{c}M_1^i,N_1^i,I_1^i\\M_2^i,N_2^i,I_2^i\\U^i,W^i,T^i,V^i\\i\in\mathbb{Z}\end{array}}\prod_{i\in\mathbb{Z}}
g^{M_1^i}_{V^iU^i}g^{N_1^i}_{T^iW^i}g^{I_2^i}_{V^iT^i}g^{M_2^i}_{U^iW^i}g^{A_1^i}_{I_1^iM_1^i}g^{A_2^i}_{N_1^iI_1^{i-1}}
g^{Y^i}_{M_2^iN_2^i}g^{A_3^i}_{N_2^i I_2^{i-1}}\\
&&\frac{a_{I_1^i}a_{N_2^i}a_{U^i}a_{V^i}a_{T^i}a_{W^i}}{a_{A_1^i}a_{A_2^i}a_{A_3^i}}\frac{\prod_{k\geq 2}\left(\langle\widehat{A_2^{i+k}},\widehat{A_1^{i}}\rangle\langle\widehat{A_3^{i+k}},\widehat{X^{i}}\rangle\right)^{(-1)^k}}
{\langle\widehat{V^i},\widehat{W^i}\rangle\langle\widehat{N_1^{i+1}},\widehat{M_1^i}\rangle\langle\widehat{N_2^{i+1}},\widehat{M_2^i}\rangle}[Y^\bullet]\\
\end{eqnarray*}
\begin{eqnarray*}
&=&\sum_{Y^\bullet\in\Iso(\cd_0(\ca))}\sum_{\begin{array}{c}I_1^i,N_2^i\\U^i,W^i,T^i,V^i\\i\in\mathbb{Z}\end{array}}\prod_{i\in\mathbb{Z}}
g^{A_1^i}_{I_1^iV^iU^i}g^{A_2^i}_{T^iW^iI_1^{i-1}}
g^{A_3^i}_{N_2^i V^{i-1}T^{i-1}}g^{Y^i}_{U^iW^iN_2^i}\\
&&\frac{a_{I_1^i}a_{N_2^i}a_{U^i}a_{V^i}a_{T^i}a_{W^i}}{a_{A_1^i}a_{A_2^i}a_{A_3^i}}\frac{\prod_{k\geq 2}\left(\langle\widehat{A_2^{i+k}},\widehat{A_1^{i}}\rangle\langle\widehat{A_3^{i+k}},\widehat{X^{i}}\rangle\right)^{(-1)^k}}
{\langle\widehat{V^i},\widehat{W^i}\rangle\langle\widehat{N_1^{i+1}},\widehat{M_1^i}\rangle\langle\widehat{N_2^{i+1}},\widehat{M_2^i}\rangle}[Y^\bullet],
\end{eqnarray*}
where $\widehat{X^i}=\widehat{U^i}+\widehat{V^i}+\widehat{T^i}+\widehat{W^i},
\widehat{N_1^{i}}=\widehat{T^i}+\widehat{W^i},
\widehat{M_1^i}=\widehat{U^i}+\widehat{V^i}$ and $\widehat{M_2^i}=\widehat{U^i}+\widehat{W^i}$ for each $i\in\mathbb{Z}$.

Similarly, we also obtain that
\begin{eqnarray*}
&&[A_1^\bullet]\left([A_2^\bullet][A_3^\bullet]\right)\\
&=&\sum_{X^\bullet\in\Iso(\cd_0(\ca))}\left(
\sum_{\begin{array}{c}P_1^i,Q_1^i,T^i\\i\in\mathbb{Z}\end{array}}
\left(\prod_{i\in\mathbb{Z}}g^{A_2^i}_{T^iP_1^i}g^{X^i}_{P_1^iQ_1^i}g^{A_3^i}_{Q_1^iT^{i-1}}
\frac{a_{T^i}a_{P_1^i}a_{Q_1^i}}{a_{A_2^i}a_{A_3^i}}\frac{\prod_{k\geq 2}\langle\widehat{A_3^{i+k}},\widehat{A_2^i}\rangle^{(-1)^k}}{\langle \widehat{Q_1^{i+1}},\widehat{P_1^i}\rangle}\right)\right)[A_1^\bullet][X^\bullet]\\
&=&\sum_{Y^\bullet,X^\bullet\in\Iso(\cd_0(\ca))}
\sum_{\begin{array}{c}P_1^i,Q_1^i,T^i\\U^i,Q_2^i,R_2^i\\i\in\mathbb{Z}\end{array}}
\prod_{i\in\mathbb{Z}}g^{A_2^i}_{T^iP_1^i}g^{X^i}_{P_1^iQ_1^i}g^{A_3^i}_{Q_1^iT^{i-1}}
g^{A_1^i}_{R_2^iU^i}g^{Y^i}_{U^iQ_2^i}g^{X^i}_{Q_2^iR_2^{i-1}}
\\
&&\frac{a_{T^i}a_{P_1^i}a_{Q_1^i}a_{R_2^i}a_{U^i}a_{Q_2^i}}{a_{A_1^i}a_{A_2^i}a_{A_3^i}a_{X^i}}\frac{\prod_{k\geq 2}\left(\langle\widehat{A_3^{i+k}},\widehat{A_2^i}\rangle\langle \widehat{X^{i+k}},\widehat{A_1^i}\rangle\right)^{(-1)^k}}{\langle \widehat{Q_1^{i+1}},\widehat{P_1^i}\rangle\langle \widehat{Q_2^{i+1}},\widehat{U^i}\rangle}[Y^\bullet]\\
&=&\sum_{Y^\bullet\in\Iso(\cd_0(\ca))}
\sum_{\begin{array}{c}P_1^i,Q_1^i,T^i\\U^i,Q_2^i,R_2^i\\N_2^i,V^i,I_1^i,W^i\\i\in\mathbb{Z}\end{array}}
\prod_{i\in\mathbb{Z}}g^{P_1^i}_{W^iI_1^{i-1}}g^{Q_1^i}_{N_2^iV^{i-1}}g^{Q_2^i}_{W^iN_2^i}g^{R_2^{i-1}}_{I_1^{i-1}V^{i-1}}
g^{A_1^i}_{R_2^iU^i}g^{A_2^i}_{T^iP_1^i}g^{A_3^i}_{Q_1^iT^{i-1}}g^{Y^i}_{U^iQ_2^i}\\
&&\frac{a_{T^i}a_{U^i}a_{W^i}a_{I_1^i}a_{N_2^i}a_{V^i}}{a_{A_1^i}a_{A_2^i}a_{A_3^i}}\frac{\prod_{k\geq 2}\left(\langle\widehat{A_3^{i+k}},\widehat{A_2^i}\rangle\langle\widehat{X^{i+k}},\widehat{A_1^i}\rangle\right)^{(-1)^k}}{\langle \widehat{Q_1^{i+1}},\widehat{P_1^i}\rangle\langle \widehat{Q_2^{i+1}},\widehat{U^i}\rangle\langle\widehat{W^i},\widehat{V^{i-1}}\rangle}[Y^\bullet]\\
&=&\sum_{Y^\bullet\in\Iso(\cd_0(\ca))}
\sum_{\begin{array}{c}T^i,U^i\\N_2^i,V^i,I_1^i,W^i\\i\in\mathbb{Z}\end{array}}
\prod_{i\in\mathbb{Z}}g^{A_1^i}_{I_1^{i}V^{i}U^i}g^{A_2^i}_{T^iW^iI_1^{i-1}}
g^{A_3^i}_{N_2^iV^{i-1}T^{i-1}}g^{Y^i}_{U^iW^iN_2^i}\\
&&\frac{a_{T^i}a_{U^i}a_{W^i}a_{I_1^i}a_{N_2^i}a_{V^i}}{a_{A_1^i}a_{A_2^i}a_{A_3^i}}\frac{\prod_{k\geq 2}\left(\langle\widehat{A_3^{i+k}},\widehat{A_2^i}\rangle\langle\widehat{X^{i+k}},\widehat{A_1^i}\rangle\right)^{(-1)^k}}{\langle \widehat{Q_1^{i+1}},\widehat{P_1^i}\rangle\langle \widehat{Q_2^{i+1}},\widehat{U^i}\rangle\langle\widehat{W^i},\widehat{V^{i-1}}\rangle}[Y^\bullet],
\end{eqnarray*}
where $\widehat{X^{i}}=\widehat{V^{i-1}}+\widehat{I_1^{i-1}}+\widehat{N_2^{i}}+\widehat{W^{i}}, \widehat{Q_1^{i}}=\widehat{V^{i-1}}+\widehat{N_2^{i}},
\widehat{P_1^i}=\widehat{I_1^{i-1}}+\widehat{W^i}$ and $\widehat{Q_2^{i}}=\widehat{N_2^{i}}+\widehat{W^i}$ for each $i\in\mathbb{Z}$.

To prove $\left([A_1^\bullet][A_2^\bullet]\right)[A_3^\bullet]=[A_1^\bullet]\left([A_2^\bullet][A_3^\bullet]\right)$, it suffices to show the following identity
\begin{eqnarray}
&&\prod_{i\in\mathbb{Z}}\frac{\prod_{k\geq 2}\left(\langle\widehat{A_2^{i+k}},\widehat{A_1^{i}}\rangle
\langle\widehat{A_3^{i+k}},\widehat{X^{i}}\rangle\right)^{(-1)^k}}
{\langle\widehat{V^i},\widehat{W^i}\rangle\langle\widehat{N_1^{i+1}},\widehat{M_1^i}\rangle
\langle\widehat{N_2^{i+1}},\widehat{M_2^i}\rangle}\label{identity of bilinear form}\\
&=&\prod_{i\in\mathbb{Z}}\frac{\prod_{k\geq 2}\left(\langle\widehat{A_3^{i+k}},\widehat{A_2^i}\rangle\langle \widehat{X^{i+k}},\widehat{A_1^i}\rangle\right)^{(-1)^k}}{\langle \widehat{V^i}+\widehat{N_2^{i+1}},\widehat{I_1^{i-1}}+\widehat{W^i}\rangle\langle \widehat{N_2^{i+1}}+\widehat{W^{i+1}},\widehat{U^i}\rangle\langle\widehat{W^i},\widehat{V^{i-1}}\rangle}\nonumber
\end{eqnarray}
for each $T^i,U^i,N_2^i,V^{i},I_1^{i},W^i\in\Iso(\ca)$ and $i\in\mathbb{Z}$.

Since $\widehat{X^{i}}=\widehat{U^i}+\widehat{V^i}+\widehat{W^i}+\widehat{T^i}=(\widehat{A_1^i}+\widehat{A_2^i})-(\widehat{I_1^{i-1}}+\widehat{I_1^{i}})$ and $\widehat{X^{i+k}}=\widehat{I_1^{i+k-1}}+\widehat{W^{i+k}}+\widehat{V^{i+k-1}}+\widehat{N_2^{i+k}}=(\widehat{A_2^{i+k}}+\widehat{A_3^{i+k}})-(\widehat{T^{i+k-1}}+\widehat{T^{i+k}})$,
then due to the following identities $$\prod_{i\in\mathbb{Z}}\prod_{k\geq 2}\langle\widehat{A_3^{i+k}},\widehat{I_1^{i-1}}+\widehat{I_1^{i}}\rangle^{(-1)^k}=
\prod_{i\in\mathbb{Z}}\langle\widehat{A_3^{i+2}},\widehat{I_1^{i}}\rangle$$
and
$$\prod_{i\in\mathbb{Z}}\prod_{k\geq 2}\langle\widehat{T^{i+k-1}}+\widehat{T^{i+k}},\widehat{A_1^i}\rangle^{(-1)^k}=\prod_{i\in\mathbb{Z}}
\langle\widehat{T^{i+1}},\widehat{A_1^i}\rangle$$ we get that
\begin{eqnarray*}
\text{LHS of}~(\ref{identity of bilinear form})&=&
\prod_{i\in\mathbb{Z}}\frac{\prod_{k\geq 2}\left(\langle\widehat{A_2^{i+k}},\widehat{A_1^{i}}\rangle
\langle\widehat{A_3^{i+k}},\widehat{A_1^i}+\widehat{A_2^i}\rangle\right)^{(-1)^k}}
{\prod_{k\geq 2}\langle\widehat{A_3^{i+k}},\widehat{I_1^{i-1}}+\widehat{I_1^{i}}\rangle^{(-1)^k}\langle\widehat{V^i},\widehat{W^i}\rangle\langle\widehat{N_1^{i+1}},\widehat{M_1^i}\rangle
\langle\widehat{N_2^{i+1}},\widehat{M_2^i}\rangle}\\
&=&\prod_{i\in\mathbb{Z}}\frac{\prod_{k\geq 2}\left(\langle\widehat{A_2^{i+k}},\widehat{A_1^{i}}\rangle
\langle\widehat{A_3^{i+k}},\widehat{A_1^i}+\widehat{A_2^i}\rangle\right)^{(-1)^k}}
{\langle\widehat{A_3^{i+2}},\widehat{I_1^{i}}\rangle\langle\widehat{V^i},\widehat{W^i}\rangle\langle\widehat{N_1^{i+1}},\widehat{M_1^i}\rangle
\langle\widehat{N_2^{i+1}},\widehat{M_2^i}\rangle}\\
&=&\prod_{i\in\mathbb{Z}}\frac{\prod_{k\geq 2}\left(\langle\widehat{A_2^{i+k}},\widehat{A_1^{i}}\rangle
\langle\widehat{A_3^{i+k}},\widehat{A_1^i}+\widehat{A_2^i}\rangle\right)^{(-1)^k}}
{\langle\widehat{A_3^{i+2}},\widehat{I_1^{i}}\rangle\langle\widehat{V^i},\widehat{W^i}\rangle
\langle\widehat{W^{i+1}}+\widehat{T^{i+1}},\widehat{U^i}+\widehat{V^i}\rangle
\langle\widehat{N_2^{i+1}},\widehat{U^i}+\widehat{W^i}\rangle}
\end{eqnarray*}
and
\begin{eqnarray*}
&&\text{RHS of}~(\ref{identity of bilinear form})\\
&=&\prod_{i\in\mathbb{Z}}\frac{\prod_{k\geq 2}\left(\langle\widehat{A_3^{i+k}},\widehat{A_2^i}\rangle\langle \widehat{A_2^{i+k}}+\widehat{A_3^{i+k}},\widehat{A_1^i}\rangle\right)^{(-1)^k}}
{\prod_{k\geq2}\langle\widehat{T^{i+k-1}}+\widehat{T^{i+k}},\widehat{A_1^i}\rangle^{(-1)^k}\langle \widehat{V^i}+\widehat{N_2^{i+1}},\widehat{I_1^{i-1}}+\widehat{W^i}\rangle\langle \widehat{N_2^{i+1}}+\widehat{W^{i+1}},\widehat{U^i}\rangle\langle\widehat{W^i},\widehat{V^{i-1}}\rangle}\\
\end{eqnarray*}
\begin{eqnarray*}
&=&\prod_{i\in\mathbb{Z}}\frac{\prod_{k\geq 2}\left(\langle\widehat{A_3^{i+k}},\widehat{A_2^i}\rangle\langle \widehat{A_2^{i+k}}+\widehat{A_3^{i+k}},\widehat{A_1^i}\rangle\right)^{(-1)^k}}
{\langle\widehat{T^{i+1}},\widehat{A_1^i}\rangle\langle\widehat{V^i}+\widehat{N_2^{i+1}},\widehat{I_1^{i-1}}+\widehat{W^i}\rangle\langle \widehat{N_2^{i+1}}+\widehat{W^{i+1}},\widehat{U^i}\rangle\langle\widehat{W^i},\widehat{V^{i-1}}\rangle}.
\end{eqnarray*}
Thus we only need to show that
\begin{eqnarray*}
&&\prod_{i\in\mathbb{Z}}\langle\widehat{A_3^{i+2}},\widehat{I_1^{i}}\rangle\langle\widehat{V^i},\widehat{W^i}\rangle
\langle\widehat{W^{i+1}}+\widehat{T^{i+1}},\widehat{U^i}+\widehat{V^i}\rangle
\langle\widehat{N_2^{i+1}},\widehat{U^i}+\widehat{W^i}\rangle\\
&=&\prod_{i\in\mathbb{Z}}\langle\widehat{T^{i+1}},\widehat{A_1^i}\rangle\langle\widehat{V^i}+\widehat{N_2^{i+1}},\widehat{I_1^{i-1}}+\widehat{W^i}\rangle\langle \widehat{N_2^{i+1}}+\widehat{W^{i+1}},\widehat{U^i}\rangle\langle\widehat{W^i},\widehat{V^{i-1}}\rangle,\end{eqnarray*}
which is equivalent to the following identity
$$\prod_{i\in\mathbb{Z}}\langle\widehat{A_3^{i+2}},\widehat{I_1^{i}}\rangle
\langle\widehat{T^{i+1}},\widehat{U^i}+\widehat{V^i}\rangle=
\prod_{i\in\mathbb{Z}}\langle\widehat{T^{i+1}},\widehat{I_1^{i}}+\widehat{U^i}+\widehat{V^i}\rangle
\langle\widehat{V^i},\widehat{I_1^{i-1}}\rangle
\langle\widehat{N_2^{i+1}},\widehat{I_1^{i-1}}\rangle.$$
And this is a direct consequence of
$\widehat{A_3^{i+2}}=\widehat{T^{i+1}}+\widehat{V^{i+1}}+\widehat{N_2^{i+2}}$.
This finishes the proof of $$\left([A_1^\bullet][A_2^\bullet]\right)[A_3^\bullet]=[A_1^\bullet]\left([A_2^\bullet][A_3^\bullet]\right)$$ in the case of $\mathcal{D}_0(\mathcal{A})$.

Assume that $t$ is an odd positive integer, let $A_1^\bullet=(A_1^i, 0)_{i\in\mathbb{Z}_t}, A_2^\bullet=(A_2^i, 0)_{i\in\mathbb{Z}_t}, A_3^\bullet=(A_3^i, 0)_{i\in\mathbb{Z}_t}\in\cd_t(\ca)$.
It suffices to show
\[\sum_{X^\bullet\in\Iso(\cd_t(\mathcal{A}))}h_{A_1^\bullet A_2^\bullet}^{X^\bullet}h_{X^\bullet A_3^\bullet}^{Y^\bullet}=\sum_{X^\bullet\in\Iso(\cd_t(\mathcal{A}))}h_{A_1^\bullet X^\bullet}^{Y^\bullet}h_{A_2^\bullet A_3^\bullet}^{X^\bullet}\]
for fixed $Y^\bullet=(Y^i, 0)_{i\in\mathbb{Z}_t}\in\Iso(\mathcal{D}_t(\mathcal{A}))$.

By definition and Theorem \ref{Green's formula} we have that
\begin{eqnarray*}
&&\sum_{X^\bullet\in\Iso(\cd_t(\mathcal{A}))}h_{A_1^\bullet A_2^\bullet}^{X^\bullet}h_{X^\bullet A_3^\bullet}^{Y^\bullet}\\
&=&\sum_{X^\bullet\in\Iso(\cd_t(\mathcal{A}))}\sum_{\begin{array}{c}
U^i,F^i,W^i,\\D^i,J^i,L^i\in\Iso(\mathcal{A})\\i\in\mathbb{Z}_t\end{array}}
\left(\prod_{i\in\mathbb{Z}_t}
\frac{g_{U^iF^i}^{A_2^i}g_{W^iU^i}^{X^i}g_{F^{i+1}W^i}^{A_1^i}g_{D^iJ^i}^{A_3^i}g_{L^iD^i}^{Y^i}g_{J^{i+1}L^i}^{X^i}}
{\langle \widehat{A_1^i},\widehat{F^i}\rangle\langle \widehat{F^{i+1}},\widehat{U^i}\rangle\langle \widehat{X^i},\widehat{J^i}\rangle\langle \widehat{J^{i+1}},\widehat{D^i}\rangle}\right.\\
&&\left.\frac{a_{U^i}a_{F^i}a_{W^i}a_{D^i}a_{J^i}a_{L^i}}{a_{X^i}a_{Y^i}}\sqrt{\langle \widehat{A_1^i},\widehat{A_2^i}\rangle\langle \widehat{X^i},\widehat{A_3^i}\rangle\prod_{k=1}^{t-1}\langle \widehat{A_1^{i+k}},\widehat{A_2^i}\rangle^{(-1)^{k+1}}\langle \widehat{X^{i+k}},\widehat{A_3^i}\rangle^{(-1)^{k+1}}}\right)\\
\end{eqnarray*}
\begin{eqnarray*}
&=&\sum_{\begin{array}{c}
U^i,F^i,W^i,\\D^i,J^i,L^i\\E
^i,R^i,S^i,G^i\in\Iso(\mathcal{A})\\i\in\mathbb{Z}_t\end{array}}\left(\prod_{i\in\mathbb{Z}_t}
\frac{g_{F^{i+1}W^i}^{A_1^i}g_{U^iF^i}^{A_2^i}g_{D^iJ^i}^{A_3^i}g_{L^iD^i}^{Y^i}g_{G^{i+1}R^i}^{W^i}g_{S^{i+1}E^i}^{U^i}g_{G^{i+1}S^{i+1}}^{J^{i+1}}
g_{R^{i}E^i}^{L^i}}
{a_{Y^i}\langle \widehat{A_1^i},\widehat{F^i}\rangle\langle \widehat{F^{i+1}},\widehat{U^i}\rangle\langle \widehat{X^i},\widehat{J^i}\rangle\langle \widehat{J^{i+1}},\widehat{D^i}\rangle\langle \widehat{G^{i+1}},\widehat{E^i}\rangle}\right.\\
&&\left.a_{D^i}a_{E^i}a_{F^i}a_{G^i}a_{R^i}a_{S^i}\sqrt{\langle \widehat{A_1^i},\widehat{A_2^i}\rangle\langle \widehat{X^i},\widehat{A_3^i}\rangle\prod_{k=1}^{t-1}\langle \widehat{A_1^{i+k}},\widehat{A_2^i}\rangle^{(-1)^{k+1}}\langle \widehat{X^{i+k}},\widehat{A_3^i}\rangle^{(-1)^{k+1}}}\right)\\
&=&\sum_{\begin{array}{c}
D^i,E^i,F^i,\\
G^i,R^i,S^i\in\Iso(\mathcal{A})\\i\in\mathbb{Z}_t\end{array}}
\left(\prod_{i\in\mathbb{Z}_t}
\frac{g_{F^{i+1}G^{i+1}R^i}^{A_1^i}g_{S^{i+1}E^iF^i}^{A_2^i}g_{D^iG^{i}S^{i}}^{A_3^i}g_{R^{i}E^iD^i}^{Y^i}a_{D^i}a_{E^i}a_{F^i}a_{G^i}a_{R^i}a_{S^i}}
{a_{Y^i}}\right.\\
&&\left.\frac{\sqrt{\langle \widehat{A_1^i},\widehat{A_2^i}\rangle\langle \widehat{X^i},\widehat{A_3^i}\rangle\prod_{k=1}^{t-1}\langle \widehat{A_1^{i+k}},\widehat{A_2^i}\rangle^{(-1)^{k+1}}\langle \widehat{X^{i+k}},\widehat{A_3^i}\rangle^{(-1)^{k+1}}}}{\langle \widehat{A_1^i},\widehat{F^i}\rangle\langle \widehat{F^{i+1}},\widehat{U^i}\rangle\langle \widehat{X^i},\widehat{J^i}\rangle\langle \widehat{J^{i+1}},\widehat{D^i}\rangle\langle \widehat{G^{i+1}},\widehat{E^i}\rangle}\right)\\
\end{eqnarray*}
and
\begin{eqnarray*}
&&\sum_{X^\bullet\in\Iso(\cd_t(\mathcal{A}))}h_{A_1^\bullet X^\bullet }^{Y^\bullet}h_{A_2^\bullet A_3^\bullet}^{X^\bullet}\\
&=&\sum_{X^\bullet\in\Iso(\cd_t(\mathcal{A}))}\sum_{\begin{array}{c}
M^i,N^i,S^i,\\P^i,Q^i,R^i\in\Iso(\mathcal{A})\\i\in\mathbb{Z}_t\end{array}}
\left(\prod_{i\in\mathbb{Z}_t}
\frac{g_{N^iS^i}^{A_3^i}g_{M^iN^i}^{X^i}g_{S^{i+1}M^i}^{A_2^i}g_{P^{i}Q^i}^{X^i}g_{R^iP^i}^{Y^i}g_{Q^{i+1}R^i}^{A_1^i}}
{\langle \widehat{A_2^i},\widehat{S^i}\rangle\langle \widehat{S^{i+1}},\widehat{N^i}\rangle
\langle \widehat{A_1^i},\widehat{Q^i}\rangle\langle \widehat{Q^{i+1}},\widehat{P^i}\rangle}\right.\\
&&\left.\frac{a_{M^i}a_{N^i}a_{S^i}a_{P^i}a_{Q^i}a_{R^i}}{a_{X^i}a_{Y^i}}\sqrt{\langle \widehat{A_2^i},\widehat{A_3^i}\rangle
\langle \widehat{A_1^i},\widehat{X^i}\rangle\prod_{k=1}^{t-1}\langle \widehat{A_2^{i+k}},\widehat{A_3^i}\rangle^{(-1)^{k+1}}
\langle\widehat{A_1^{i+k}}, \widehat{X^i}\rangle^{(-1)^{k+1}}}\right)\\
&=&\sum_{\begin{array}{c}
M^i,N^i,S^i,\\P^i,Q^i,R^i\\
D^i,E^i,F^i,G^i\in\Iso(\mathcal{A})\\i\in\mathbb{Z}_t\end{array}}
\left(\prod_{i\in\mathbb{Z}_t}
\frac{g_{Q^{i+1}R^i}^{A_1^i}g_{S^{i+1}M^i}^{A_2^i}g_{N^iS^i}^{A_3^i}g_{R^iP^i}^{Y^i}g_{F^iG^i}^{Q^i}g_{E^{i}F^i}^{M^i}
g_{D^{i}G^i}^{N^i}g_{E^{i}D^i}^{P^i}}
{a_{Y^i}\langle \widehat{A_2^i},\widehat{S^i}\rangle\langle \widehat{S^{i+1}},\widehat{N^i}\rangle
\langle \widehat{A_1^i},\widehat{Q^i}\rangle\langle \widehat{Q^{i+1}},\widehat{P^i}\rangle\langle\widehat{E^i},\widehat{G^i}\rangle}\right.\\
&&\left.a_{D^i}a_{E^i}a_{F^i}a_{G^i}a_{R^i}a_{S^i}\sqrt{\langle \widehat{A_2^i},\widehat{A_3^i}\rangle
\langle \widehat{A_1^i},\widehat{X^i}\rangle\prod_{k=1}^{t-1}\langle \widehat{A_2^{i+k}},\widehat{A_3^i}\rangle^{(-1)^{k+1}}
\langle\widehat{A_1^{i+k}}, \widehat{X^i}\rangle^{(-1)^{k+1}}}\right)\\
&=&\sum_{\begin{array}{c}
D^i,E^i,F^i,\\G^i,R^i,S^i\in\Iso(\mathcal{A})\\i\in\mathbb{Z}_t\end{array}}
\left(\prod_{i\in\mathbb{Z}_t}
\frac{g_{F^{i+1}G^{i+1}R^i}^{A_1^i}g_{S^{i+1}E^iF^i}^{A_2^i}g_{D^iG^{i}S^{i}}^{A_3^i}g_{R^{i}E^iD^i}^{Y^i}a_{D^i}a_{E^i}a_{F^i}a_{G^i}a_{R^i}a_{S^i}}
{a_{Y^i}}\right.\\
\end{eqnarray*}
\begin{eqnarray*}
&&\left.\frac{\sqrt{\langle \widehat{A_2^i},\widehat{A_3^i}\rangle
\langle \widehat{A_1^i},\widehat{X^i}\rangle\prod_{k=1}^{t-1}\langle \widehat{A_2^{i+k}},\widehat{A_3^i}\rangle^{(-1)^{k+1}}
\langle\widehat{A_1^{i+k}}, \widehat{X^i}\rangle^{(-1)^{k+1}}}}{\langle \widehat{A_2^i},\widehat{S^i}\rangle\langle \widehat{S^{i+1}},\widehat{N^i}\rangle
\langle \widehat{A_1^i},\widehat{Q^i}\rangle\langle \widehat{Q^{i+1}},\widehat{P^i}\rangle\langle\widehat{E^i},\widehat{G^i}\rangle}\right).
\end{eqnarray*}
So it suffices to show
\begin{eqnarray*}
&&\prod_{i\in\mathbb{Z}_t}\frac{\sqrt{\langle \widehat{A_1^i},\widehat{A_2^i}\rangle\langle \widehat{X^i},\widehat{A_3^i}\rangle\prod_{k=1}^{t-1}\langle \widehat{A_1^{i+k}},\widehat{A_2^i}\rangle^{(-1)^{k+1}}\langle \widehat{X^{i+k}},\widehat{A_3^i}\rangle^{(-1)^{k+1}}}}{\langle \widehat{A_1^i},\widehat{F^i}\rangle\langle \widehat{F^{i+1}},\widehat{U^i}\rangle\langle \widehat{X^i},\widehat{J^i}\rangle\langle \widehat{J^{i+1}},\widehat{D^i}\rangle\langle \widehat{G^{i+1}},\widehat{E^i}\rangle}\\
&=&\prod_{i\in\mathbb{Z}_t}\frac{\sqrt{\langle \widehat{A_2^i},\widehat{A_3^i}\rangle
\langle \widehat{A_1^i},\widehat{X^i}\rangle\prod_{k=1}^{t-1}\langle \widehat{A_2^{i+k}},\widehat{A_3^i}\rangle^{(-1)^{k+1}}
\langle\widehat{A_1^{i+k}}, \widehat{X^i}\rangle^{(-1)^{k+1}}}}{\langle \widehat{A_2^i},\widehat{S^i}\rangle\langle \widehat{S^{i+1}},\widehat{N^i}\rangle
\langle \widehat{A_1^i},\widehat{Q^i}\rangle\langle \widehat{Q^{i+1}},\widehat{P^i}\rangle\langle\widehat{E^i},\widehat{G^i}\rangle}
\end{eqnarray*}
for fixed $D^i,E^i,F^i,G^i,R^i,S^i\in\Iso(\mathcal{A})$, which is equivalent to the following identity
\begin{eqnarray}
&&\prod_{i\in\mathbb{Z}_t}\frac{\langle \widehat{A_1^i},\widehat{F^i}\rangle\langle \widehat{F^{i+1}},\widehat{U^i}\rangle\langle \widehat{X^i},\widehat{J^i}\rangle\langle \widehat{J^{i+1}},\widehat{D^i}\rangle\langle \widehat{G^{i+1}},\widehat{E^i}\rangle}{\langle \widehat{A_2^i},\widehat{S^i}\rangle\langle \widehat{S^{i+1}},\widehat{N^i}\rangle
\langle \widehat{A_1^i},\widehat{Q^i}\rangle\langle \widehat{Q^{i+1}},\widehat{P^i}\rangle\langle\widehat{E^i},\widehat{G^i}\rangle}\label{identity of associativity}\\
&=&\prod_{i\in\mathbb{Z}_t}\sqrt{\frac{\langle \widehat{A_1^i},\widehat{A_2^i}\rangle\langle \widehat{X^i},\widehat{A_3^i}\rangle\prod_{k=1}^{t-1}\langle \widehat{A_2^{i+k}},\widehat{A_3^i}\rangle^{(-1)^{k}}
\langle\widehat{A_1^{i+k}}, \widehat{X^i}\rangle^{(-1)^{k}}}
{\langle \widehat{A_2^i},\widehat{A_3^i}\rangle
\langle \widehat{A_1^i},\widehat{X^i}\rangle\prod_{k=1}^{t-1}\langle \widehat{A_1^{i+k}},\widehat{A_2^i}\rangle^{(-1)^{k}}\langle \widehat{X^{i+k}},\widehat{A_3^i}\rangle^{(-1)^{k}}}}\nonumber.
\end{eqnarray}
And we easily get that

\begin{eqnarray*}
&&\text{LHS of}~(\ref{identity of associativity})\\
&=&\prod_{i\in\mathbb{Z}_t}\frac{\langle \widehat{A_1^i},\widehat{F^i}\rangle\langle \widehat{F^{i+1}},\widehat{E^i}+\widehat{S^{i+1}}\rangle\langle \widehat{E^i}+\widehat{R^i}+\widehat{S^{i+1}}+\widehat{G^{i+1}},\widehat{S^{i}}+\widehat{G^{i}}\rangle\langle \widehat{S^{i+1}}+\widehat{G^{i+1}},\widehat{D^i}\rangle\langle \widehat{G^{i+1}},\widehat{E^i}\rangle}
{\langle\widehat{E^i},\widehat{G^i}\rangle\langle \widehat{A_2^i},\widehat{S^i}\rangle\langle \widehat{S^{i+1}},\widehat{D^i}+\widehat{G^{i}}\rangle
\langle \widehat{A_1^i},\widehat{F^i}+\widehat{G^{i}}\rangle\langle \widehat{F^{i+1}}+\widehat{G^{i+1}},\widehat{D^i}+\widehat{E^i}\rangle}\\
&=&\prod_{i\in\mathbb{Z}_t}\frac{\langle \widehat{F^{i+1}},\widehat{E^i}+\widehat{S^{i+1}}\rangle\langle \widehat{E^i}+\widehat{R^i}+\widehat{S^{i+1}}+\widehat{G^{i+1}},\widehat{S^{i}}+\widehat{G^{i}}\rangle\langle \widehat{S^{i+1}},\widehat{D^i}\rangle\langle \widehat{G^{i+1}},\widehat{D^i}+\widehat{E^i}\rangle}
{\langle\widehat{E^i}+\widehat{S^{i+1}}+\widehat{A_1^i},\widehat{G^i}\rangle\langle \widehat{A_2^i},\widehat{S^i}\rangle\langle \widehat{S^{i+1}},\widehat{D^i}\rangle
\langle \widehat{F^{i+1}}+\widehat{G^{i+1}},\widehat{D^i}+\widehat{E^i}\rangle}\\
&=&\prod_{i\in\mathbb{Z}_t}
\frac{\langle\widehat{R^i}+\widehat{G^{i+1}},\widehat{S^i}\rangle\langle\widehat{F^{i+1}},\widehat{E^i}+\widehat{S^{i+1}}\rangle}
{\langle \widehat{F^{i+1}},\widehat{G^i}\rangle
\langle\widehat{F^{i}},\widehat{S^i}\rangle
\langle \widehat{F^{i+1}},\widehat{D^i}+\widehat{E^i}\rangle}\\
&=&\prod_{i\in\mathbb{Z}_t}
\frac{\langle\widehat{R^i}+\widehat{G^{i+1}},\widehat{S^i}\rangle}
{\langle\widehat{F^{i+1}},\widehat{G^i}+\widehat{D^i}\rangle}.
\end{eqnarray*}

To calculate the RHS of (\ref{identity of associativity}), we should note the following facts
\begin{eqnarray*}
\prod_{i\in\mathbb{Z}_t}\frac{\langle \widehat{A_1^i},\widehat{A_2^i}\rangle\langle \widehat{X^i},\widehat{A_3^i}\rangle}{\langle \widehat{A_2^i},\widehat{A_3^i}\rangle
\langle \widehat{A_1^i},\widehat{X^i}\rangle}&=&\prod_{i\in\mathbb{Z}_t}\frac{\langle \widehat{A_1^i},\widehat{A_2^i}\rangle\langle \widehat{E^i}+\widehat{R^i}+\widehat{S^{i+1}}+\widehat{G^{i+1}},\widehat{A_3^i}\rangle}
{\langle \widehat{A_2^i},\widehat{A_3^i}\rangle
\langle\widehat{A_1^i},\widehat{D^i}+\widehat{E^i}+\widehat{F^{i}}+\widehat{G^{i}}\rangle}\\
&=&\prod_{i\in\mathbb{Z}_t}\frac{\langle \widehat{A_1^i},\widehat{S^{i+1}}\rangle\langle\widehat{R^i}+\widehat{G^{i+1}},\widehat{A_3^i}\rangle}
{\langle \widehat{A_1^i},\widehat{D^{i}}+\widehat{G^{i}}\rangle\langle \widehat{F^i},\widehat{A_3^i}\rangle}\\
&=&\prod_{i\in\mathbb{Z}_t}
\frac{\langle\widehat{F^{i+1}}+\widehat{G^{i+1}}+\widehat{R^i},\widehat{S^{i+1}}\rangle
\langle\widehat{R^i}+\widehat{G^{i+1}},\widehat{D^i}+\widehat{G^i}+\widehat{S^i}\rangle}
{\langle\widehat{F^{i+1}}+\widehat{G^{i+1}}+\widehat{R^i},\widehat{D^{i}}+\widehat{G^{i}}\rangle\langle \widehat{F^i},\widehat{D^i}+\widehat{G^i}+\widehat{S^i}\rangle}\\
&=&\prod_{i\in\mathbb{Z}_t}
\frac{\langle\widehat{R^i}+\widehat{G^{i+1}},\widehat{S^i}+\widehat{S^{i+1}}\rangle}
{\langle\widehat{F^{i+1}}+\widehat{F^i},\widehat{D^{i}}+\widehat{G^{i}}\rangle}\\
\end{eqnarray*}
and
\begin{eqnarray*}
&&\prod_{i\in\mathbb{Z}_t}\prod_{k=1}^{t-1}\left(\frac{\langle \widehat{A_2^{i+k}},\widehat{A_3^i}\rangle
\langle\widehat{A_1^{i+k}}, \widehat{X^i}\rangle}
{\langle \widehat{A_1^{i+k}},\widehat{A_2^i}\rangle\langle \widehat{X^{i+k}},\widehat{A_3^i}\rangle}
\right)^{(-1)^k}\\
&=&\prod_{i\in\mathbb{Z}_t}\prod_{k=1}^{t-1}\left(
\frac{\langle \widehat{F^{i+k}},\widehat{A_3^i}\rangle
\langle\widehat{A_1^{i+k}}, \widehat{D^i}+\widehat{G^i}\rangle}
{\langle \widehat{R^{i+k}}+\widehat{G^{i+k+1}},\widehat{A_3^i}\rangle\langle \widehat{A_1^{i+k}},\widehat{S^{i+1}}\rangle}
\right)^{(-1)^k}\\
&=&\prod_{i\in\mathbb{Z}_t}\prod_{k=1}^{t-1}
\left(
\frac{\langle\widehat{F^{i+k}},\widehat{D^i}+\widehat{G^i}+\widehat{S^{i}}\rangle
\langle\widehat{F^{i+k+1}}+\widehat{G^{i+k+1}}+\widehat{R^{i+k}}, \widehat{D^i}+\widehat{G^i}\rangle}
{\langle \widehat{R^{i+k}}+\widehat{G^{i+k+1}},\widehat{D^i}+\widehat{G^i}+\widehat{S^{i}}\rangle\langle \widehat{F^{i+k+1}}+\widehat{G^{i+k+1}}+\widehat{R^{i+k}},\widehat{S^{i+1}}\rangle}
\right)^{(-1)^k}\\
&=&\prod_{i\in\mathbb{Z}_t}\prod_{k=1}^{t-1}
\left(
\frac{\langle\widehat{F^{i+k}}+\widehat{F^{i+k+1}},\widehat{D^i}+\widehat{G^i}\rangle
\langle\widehat{F^{i+k}},\widehat{S^{i}}\rangle}
{\langle \widehat{R^{i+k}}+\widehat{G^{i+k+1}},\widehat{S^i}\rangle\langle \widehat{F^{i+k+1}}+\widehat{G^{i+k+1}}+\widehat{R^{i+k}},\widehat{S^{i+1}}\rangle}
\right)^{(-1)^k}\\
&=&\prod_{i\in\mathbb{Z}_t}\prod_{k=1}^{t-1}
\left(\frac{\langle\widehat{F^{i+k}}+\widehat{F^{i+k+1}},\widehat{D^i}+\widehat{G^i}\rangle}
{\widehat{R^{i+k}}+\widehat{G^{i+k+1}},\widehat{S^i}\rangle\langle \widehat{G^{i+k+1}}+\widehat{R^{i+k}},\widehat{S^{i+1}}\rangle}
\right)^{(-1)^k}
\left(
\frac{\langle\widehat{F^{i+k}},\widehat{S^{i}}\rangle}{\langle\widehat{F^{i+k+1}},\widehat{S^{i+1}}\rangle}
\right)^{(-1)^k}\\
&=&\prod_{i\in\mathbb{Z}_t}\prod_{k=1}^{t-1}
\left(\frac{\langle\widehat{F^{i+k}}+\widehat{F^{i+k+1}},\widehat{D^i}+\widehat{G^i}\rangle}
{\widehat{R^{i+k}}+\widehat{G^{i+k+1}},\widehat{S^i}\rangle\langle \widehat{R^{i+k-1}}+\widehat{G^{i+k}},\widehat{S^{i}}\rangle}
\right)^{(-1)^k}
\end{eqnarray*}
since $\prod_{i\in\mathbb{Z}_t}\prod_{k=1}^{t-1}
\left(\frac{\langle\widehat{F^{i+k}},\widehat{S^{i}}\rangle}{\langle\widehat{F^{i+k+1}},\widehat{S^{i+1}}\rangle}
\right)^{(-1)^k}=1$. We have
\begin{eqnarray*}
&&\prod_{k=1}^{t-1}\langle\widehat{F^{i+k}}+\widehat{F^{i+k+1}},\widehat{D^i}+\widehat{G^i}\rangle^{(-1)^k}\\
&=&\langle(\widehat{F^{i+2}}+\widehat{F^{i+3}})-(\widehat{F^{i+1}}+\widehat{F^{i+2}})
%+(\widehat{F^{i+4}}+\widehat{F^{i+5}})-(\widehat{F^{i+3}}+\widehat{F^{i+4}})
+\cdots+
(\widehat{F^{i+t-1}}+\widehat{F^{i+t}})-(\widehat{F^{i+t-2}}+\widehat{F^{i+t-1}}),
\widehat{D^i}+\widehat{G^i}\rangle\\
&=&\langle\widehat{F^{i}}-\widehat{F^{i+1}},\widehat{D^i}+\widehat{G^i}\rangle,
\end{eqnarray*}
similarly one can obtain that
$$\prod_{k=1}^{t-1}\langle\widehat{G^{i+k}}+\widehat{G^{i+k+1}},\widehat{S^i}\rangle^{(-1)^k}
=\langle\widehat{G^{i}}-\widehat{G^{i+1}},\widehat{S^i}\rangle$$
and
$$\prod_{k=1}^{t-1}\langle\widehat{R^{i+k-1}}+\widehat{R^{i+k}},\widehat{S^i}\rangle^{(-1)^k}
=\langle\widehat{R^{i-1}}-\widehat{R^{i}},\widehat{S^i}\rangle.$$

Thus we have that
\begin{eqnarray*}
\text{RHS of (\ref{identity of associativity})}
&=&\prod_{i\in\mathbb{Z}_t}\sqrt{\frac{\langle\widehat{R^i}+\widehat{G^{i+1}},\widehat{S^i}+\widehat{S^{i+1}}\rangle}
{\langle\widehat{F^{i+1}}+\widehat{F^i},\widehat{D^{i}}+\widehat{G^{i}}\rangle}
\frac{\langle\widehat{F^{i}}-\widehat{F^{i+1}},\widehat{D^i}+\widehat{G^i}\rangle}{\langle\widehat{G^{i}}-\widehat{G^{i+1}}+\widehat{R^{i-1}}-\widehat{R^{i}},\widehat{S^i}\rangle}}\\
&=&\prod_{i\in\mathbb{Z}_t}\sqrt{\frac{\langle\widehat{R^i}+\widehat{G^{i+1}},\widehat{S^i}+\widehat{S^{i+1}}\rangle}
{\langle\widehat{F^{i+1}}+\widehat{F^i},\widehat{D^{i}}+\widehat{G^{i}}\rangle
}
\frac{\langle\widehat{F^{i}}-\widehat{F^{i+1}},\widehat{D^i}+\widehat{G^i}\rangle
\langle\widehat{R^{i}}+\widehat{G^{i+1}},\widehat{S^i}\rangle}{\langle\widehat{G^{i}}+\widehat{R^{i-1}},\widehat{S^i}\rangle}}\\
&=&\prod_{i\in\mathbb{Z}_t}\frac{\langle\widehat{R^i}+\widehat{G^{i+1}},\widehat{S^i}\rangle}
{\langle\widehat{F^{i+1}},\widehat{G^i}+\widehat{D^i}\rangle}\\
&=&\text{LHS of (\ref{identity of associativity})}.
\end{eqnarray*}
This finishes the proof.
\end{proof}

\section{Derived Hall algebras}\label{Sec derived Hall algebras}
Let $\ct$ be a $k$-additive triangulated category. For any objects $X, Y, L\in\ct$, we denote by $\Hom_{\ct}(Y, L)_{X}$ the subset of $\Hom_{\ct}(Y, L)$ consisting of morphisms $l: Y\rightarrow L$ whose cone $\Cone(l)$ is isomorphic to $X$.
\subsection{Hall algebras associated to triangulated categories without periodicities}Let $\ct$ be a $k$-additive triangulated category with the translation $[1]$ satisfying
\begin{itemize}
\item[(i)] $\dim_k\Hom_{\ct}(X, Y)<\infty$ for any two objects $X$ and $Y$;
\item[(ii)] $\End_{\ct}(X)$ is local for any indecomposable object $X$;
\item[(iii)] $\ct$ is (left) locally finite; that is, $\sum_{i\geq 0}\dim_k\Hom_{\ct}(X[i], Y)<\infty$ for any $X$ and $Y$.
\end{itemize}

Note that the first two conditions imply the validity of the Krull-Schmidt theorem in ${\ct}$, which means that any object in ${\ct}$ can be uniquely decomposed into the direct sum of finitely many indecomposable objects up to isomorphism.

For any objects $X, Y\in\ct$, set $$\{X, Y\}:=\prod_{i>0}|\Hom_{\ct}(X[i], Y)|^{(-1)^i}.$$
Similarly, for any $X, Y, L\in\ct$, we have
$$\frac{|\Hom_{\ct}(L,X)_{Y[1]}|}{a_L}\cdot\frac{\{L, X\}}{\{L, L\}}=\frac{|\Hom_{\ct}(X, Y[1])_{L[1]}|}{a_Y}\cdot\frac{\{X,Y[1]\}}{\{Y, Y\}}$$ and
$$\frac{|\Hom_{\ct}(Y, L)_X|}{a_L}\cdot\frac{\{Y, L\}}{\{L, L\}}=\frac{|\Hom_{\ct}(X[-1], Y)_L|}{a_X}\cdot\frac{\{X[-1], Y\}}{\{X, X\}}.$$
Thus
$$\frac{|\Hom_{\ct}(L,X)_{Y[1]}|}{a_X}\cdot\frac{\{L, X\}}{\{X, X\}}=\frac{|\Hom_{\ct}(Y, L)_X|}{a_Y}\cdot\frac{\{Y, L\}}{\{Y, Y\}}$$
for any $X, Y, L\in\ct$, denote it by $F_{XY}^L$. And then we have the derived Riedtmann-Peng formula, i.e.,
$$F_{XY}^L=|\Hom_{\ct}(X,Y[1])_{L[1]}|\{X,Y[1]\}\frac{a_L\{L,L\}}{a_X\{X,X\}a_Y\{Y,Y\}}.$$

The derived Hall algebra $\cd\ch(\ct)$ of the triangulated category $\ct$ is the $\mathbb{Q}$-space with the basis $\{[X]|X\in\ct\}$ and the multiplication is defined by
$$[X][Y]=\sum_{[L]}F_{XY}^{L}[L].$$

Obviously, the bounded derived category $\cd_0(\ca)$ of hereditary abelian category $\ca$ satisfies the above conditions (i)-(iii), and we denote by $\cd\ch_0(\ca)$ the derived Hall algebra of $\cd_0(\ca)$.

In the following, for any object $M,B,A,N\in\ca$, denote by $$\gamma_{AB}^{MN}:=\frac{|V(M,B,A,N)|}{a_Aa_B},$$ where $V(M,B,A,N)$ is the set of all exact sequences
\[0\rightarrow M\rightarrow B\rightarrow A\rightarrow N\rightarrow 0,\]
also we have $\gamma_{AB}^{MN}=\sum\limits_{I\in \Iso(\ca)}g_{IM}^Bg_{NI}^A\frac{a_Ma_Na_I}{a_Aa_B}$.
The derived Hall algebra of $\cd_0(\ca)$ can be described by generators and relations as follows.
\begin{proposition}[\cite{T06}, Proposition 7.1]\label{proposition derived hall algebra}
$\cd\ch_0(\ca)$ is isomorphic to an associative and unital $\mathbb{Q}$-algebra generated by the set $$\{Z_A^{[n]}|A\in\Iso(\ca), n\in\mathbb{Z}\},$$  with the
following defining relations (\ref{relation 4.3})-(\ref{relation 4.5}).
\begin{eqnarray}
Z_A^{[n]}Z_B^{[n]}&=&\sum\limits_{C\in \Iso(\ca)}g_{AB}^CZ_C^{[n]},\label{relation 4.3}\\
Z_B^{[n]}Z_A^{[n+1]}&=&\sum\limits_{M, N\in\Iso(\ca)}\gamma_{AB}^{MN}\frac{1}{\langle \widehat{N}, \widehat{M}\rangle}Z_N^{[n+1]}Z_M^{[n]},\label{relation 4.4}\\
Z_B^{[n]}Z_A^{[m]}&=&\langle \widehat{A}, \widehat{B}\rangle^{(-1)^{m-n}}Z_A^{[m]}Z_B^{[n]}\ \text{for~} m>n+1,\label{relation 4.5}
\end{eqnarray}
for any $A, B\in\Iso(\ca)$ and $n\in\mathbb{Z}$.
\end{proposition}
By definition, one can easily get the following lemma.
\begin{lemma}\label{lemma multiplication of large item product small item}
For any $A, B\in\Iso(\ca)$ and $i>j\in\mathbb{Z}$, we have $$Z_A^{[i]}Z_B^{[j]}=Z_A^{[i]}\oplus Z_B^{[j]}$$
in $\cd\ch_0(\ca)$.
\end{lemma}

\subsection{Hall algebras associated to odd periodic triangulated categories}Let $\ct$ be a $k$-additive triangulated category with the translation $[1]$ satisfying the following conditions:
\begin{itemize}
\item[(i)] $\dim_k\Hom_{\ct}(X, Y)<\infty$ for any two objects $X$ and $Y$;
\item[(ii)] $\End_{\ct}(X)$ is local for any indecomposable object $X$;
\item[(iii)] $[1]^t=[t]\cong 1_{\ct}$ for some positive integer $t$.
\end{itemize}
The category $\ct$ is then called a $t$-periodic triangulated category and $[1]$ is called a $t$-periodic translation (or shift) functor.

In particular, Xu and Chen defined the Hall algebra associated to any odd periodic triangulated category satisfying some finiteness conditions which generalized the results in \cite{T06, XX08, XX15}. For any $X, Y, L\in\ct$, we have an analogue of To\"{e}n's formula as follows
$$\frac{|\Hom_{\ct}(L,X)_{Y[1]}|}{a_X}\left(\prod_{i=1}^{t}\frac{|\Hom_{\ct}(L[i], X)|^{(-1)^i}}{|\Hom_{\ct}(X[i], X)|^{(-1)^i}}\right)^{\frac{1}{2}}
=\frac{|\Hom_{\ct}(Y,L)_{X}|}{a_Y}\left(\prod_{i=1}^{t}\frac{|\Hom_{\ct}(Y[i], L)|^{(-1)^i}}{|\Hom_{\ct}(Y[i], Y)|^{(-1)^i}}\right)^{\frac{1}{2}},$$
and we also denote it by $F_{XY}^L$.
By abuse of notation, we shall denote by $\{X,Y\}$ the alternate multiplication $$\prod_{i=1}^{t}|\Hom_{\ct}(X[i], Y)|^{(-1)^i}.$$
Furthermore, one can get the derived Riedtmann-Peng formula as follows
$$F_{XY}^L=|\Hom_{\ct}(X,Y[1])_{L[1]}|\{X,Y[1]\}^{\frac{1}{2}}\frac{a_L\{L,L\}^{\frac{1}{2}}}{a_X\{X,X\}^{\frac{1}{2}}a_Y\{Y,Y\}^{\frac{1}{2}}},$$
see \cite{XX15, SCX18} for more details.

Set $v=\sqrt{q}$, then $F_{XY}^L\in\mathbb{Q}(v)$. Here $\mathbb{Q}(v)$ be the rational field of $\mathbb{Q}[v]$. The derived Hall algebra $\cd\ch_t(\ct)$ of an odd periodic triangulated category $\ct$ introduced by Xu and Chen in \cite{XC13} is the $\mathbb{Q}(v)$-space with the basis $\{[X]|X\in\ct\}$ and the multiplication is defined to be
$$[X][Y]=\sum_{[L]\in\Iso(\ct)}F_{XY}^L[L].$$
Obviously we can define the derived Hall algebra of $\cd_t(\ca)$ for any odd integer $t$, and denote it by $\cd\ch_t(\ca)$.
Recently, applying the derived Hall numbers of the bounded derived category $\cd_0(\ca)$, H. Zhang \cite{Zhang23} gave an explicit characterization for the multiplication structure of the odd periodic derived Hall algebra of $\ca$.

Similar to the work of Guo and Peng  in \cite{GP97}, one can find a universal PBW basis of $\cd\ch_t(\ca)$ and then describe the derived Hall algebra $\cd\ch_t(\ca)$ in terms of generators and relations as follows.

\begin{proposition}\label{proposition generators and relations 1-derived hall algebra}
$\cd\ch_1(\ca)$ is isomorphic to an associative and unital $\mathbb{Q}(v)$-algebra generated by the set $$\{Z_A|A\in\Iso(\ca)\},$$  with the defining relation
\begin{eqnarray}
Z_AZ_B&=&\sum\limits_{C\in\Iso(\ca)}\frac{|\Hom_{\cd_1(\ca)}(Z_A,Z_B)_{Z_C}|}{\sqrt{|\Hom_{\ca}(A,B)||\Ext^1_{\ca}(A,B)|}}\frac{a'_{C}}{a'_{A}a'_{B}}Z_C,\label{re1 in DH1}
\end{eqnarray}
where $a'_{X}=\sqrt{|\Hom_{\ca}(X,X)||\Ext^1_{\ca}(X,X)|}$ for any $X\in\Iso(\ca)$ and $A, B\in\Iso(\ca)$.
\end{proposition}

\begin{proof}
For any $X,Y\in\Iso(\ca)$, we have that
\begin{eqnarray*}
\{Z_X,Z_Y\}&=&\frac{1}{|\Hom_{\cd_1(\ca)}(Z_X[1],Z_Y)|}\\
&=&\frac{1}{|\bigoplus_{k\in\mathbb{Z}}\Hom_{\cd^b(\ca)}(X[1],Y[k])|}\\
&=&\frac{1}{|\Hom_{\mathcal{A}}(X,Y)||\Ext^1_{\mathcal{A}}(X,Y)|}
\end{eqnarray*}
and
\begin{eqnarray*}
a_{Z_{X}}\{Z_X,Z_X\}^{\frac{1}{2}}&=&|\Hom_{\cd_1(\ca)}(Z_X,Z_X)|\frac{1}{\sqrt{|\Hom_{\cd_1(\ca)}(Z_X[1],Z_X)|}}\\
&=&\sqrt{|\Hom_{\cd_1(\ca)}(Z_X,Z_X)|}\\
&=&\sqrt{|\Hom_{\mathcal{A}}(X,X)||\Ext^1_{\mathcal{A}}(X,X)|}\\
&=&a_{X}'.
\end{eqnarray*}
in $\cd\ch_1(\ca)$.

By Proposition \ref{proposition t-periodic iso to the homology} we know that
\[Z_AZ_B=\sum\limits_{C\in\Iso(\ca)}F_{Z_AZ_B}^{Z_C}Z_C,\]
where
\begin{eqnarray*}
F_{Z_AZ_B}^{Z_C}&=&|\Hom_{\mathcal{D}_1(\mathcal{A})}(Z_A,Z_B[1])_{Z_C[1]}|\{Z_A,Z_B[1]\}^{\frac{1}{2}}
\frac{a_{Z_C}\{Z_C,Z_C\}^{\frac{1}{2}}}{a_{Z_A}\{Z_A,Z_A\}^{\frac{1}{2}}a_{Z_B}\{Z_B,Z_B\}^{\frac{1}{2}}}\\
&=&\frac{|\Hom_{\cd_1(\ca)}(Z_A,Z_B)_{Z_C}|}{\sqrt{|\Hom_{\ca}(A,B)||\Ext^1_{\ca}(A,B)|}}\frac{a'_{C}}{a'_{A}a'_{B}}
\end{eqnarray*}
This finishes the proof.
\end{proof}

\begin{proposition}\label{proposition generators and relations 3-derived hall algebra}
$\cd\ch_3(\ca)$ is isomorphic to an associative and unital $\mathbb{Q}(v)$-algebra generated by the set $$\{Z_A^{[i]}|A\in\Iso(\ca),i\in\mathbb{Z}_3\},$$  with the defining relations as follows.
\begin{eqnarray}
Z_A^{[i]}Z_B^{[i]}&=&\sum\limits_{C\in \Iso(\ca)}g_{AB}^C\frac{1}{\sqrt{\langle \widehat{B},\widehat{A}\rangle}}Z_C^{[i]},\label{prop relation 1 in DH3}\\
Z_B^{[i]}Z_A^{[i+1]}&=&\sum\limits_{M, N\in\Iso(\ca)}\gamma_{AB}^{MN}\sqrt{\frac{\langle \widehat{A},\widehat{A}\rangle\langle\widehat{B},\widehat{B}\rangle}{\langle\widehat{M},\widehat{M}\rangle\langle\widehat{N},\widehat{N}\rangle}}\frac{1}{\sqrt{\langle\widehat{B}, \widehat{A}\rangle\langle \widehat{N}, \widehat{M}\rangle}}Z_N^{[i+1]}Z_M^{[i]},\label{prop relation 2 in DH3}
\end{eqnarray}
where $A, B\in\Iso(\ca)$ and $i\in\mathbb{Z}_3$.
\end{proposition}

\begin{proof}
We only show (\ref{prop relation 1 in DH3}), the proof of (\ref{prop relation 2 in DH3}) is similar to that in \cite{T06} and \cite{XC13}. Following \cite{T06,XX15,XC13}, one can get that
\begin{eqnarray*}
Z_A^{[i]}Z_B^{[i]}&=&\sum\limits_{C\in \Iso(\ca)}\frac{|\Ext^1_{\mathcal{A}}(A,B)_C|}{\sqrt{|\Hom_{\mathcal{A}}(A,B)||\Ext^1_{\mathcal{A}}(A,B)|}}\frac{a'_{C}}{a'_{A}a'_{B}}Z_C^{[i]}\\
&=&\sum\limits_{C\in \Iso(\ca)}\frac{g_{AB}^C|\Hom_{\mathcal{A}}(A,B)|}
{\sqrt{|\Hom_{\mathcal{A}}(A,B)||\Ext^1_{\mathcal{A}}(A,B)|}}
\frac{a_{A}a_{B}}{a_{C}}\frac{a'_{C}}{a'_{A}a'_{B}}Z_C^{[i]}\\
&=&\sum\limits_{C\in \Iso(\ca)}g_{AB}^C\frac{\sqrt{\langle \widehat{A},\widehat{B}\rangle\langle \widehat{A},\widehat{A}\rangle\langle \widehat{B},\widehat{B}\rangle}}{\sqrt{\langle \widehat{C},\widehat{C}\rangle}}Z_C^{[i]}\\
&=&\sum\limits_{C\in \Iso(\ca)}g_{AB}^C\frac{1}{\sqrt{\langle \widehat{B},\widehat{A}\rangle}}Z_C^{[i]}.
\end{eqnarray*}
This completes the proof.
\end{proof}

Similarly one can get the generators and relations of $\cd\ch_t(\ca)$ as follows.
\begin{proposition}\label{proposition generators and relations t-derived hall algebra}
$\cd\ch_t(\ca)(t>3)$ is isomorphic to an associative and unital $\mathbb{Q}(v)$-algebra generated by the set $$\{Z_A^{[i]}|A\in\Iso(\ca),i\in\mathbb{Z}_t\},$$  with the defining relations (\ref{prop relation 1 in DH3})-(\ref{relation 3 in DHt}).
\begin{eqnarray}
Z_A^{[i]}Z_B^{[j]}&=&\sqrt{(\widehat{A}, \widehat{B})^{(-1)^{j-i}}}Z_B^{[j]}Z_A^{[i]}\ \text{for any~} j=i+2,\cdots,t-1,\label{relation 3 in DHt}
\end{eqnarray}
where $A, B\in\Iso(\ca)$ and $i,j\in\mathbb{Z}_t$.
\end{proposition}

\section{Isomorphisms between $\cl_t(\ca)$ and derived Hall algebras}\label{Sec derived Hall numbers}
In this section we will work out the structure constants of complexes with zero differentials in derived Hall algebras and then show the isomorphisms between $\cl_t(\ca)$ and derived Hall algebras.
\subsection{Derived Hall numbers for bounded complexes}

For any nonzero complex $M^\bullet\in\cc_0(\ca)$,  $$M^\bullet=\cdots\rightarrow 0\rightarrow M^l\rightarrow\cdots\rightarrow M^r\rightarrow 0\rightarrow\cdots,$$
where $M^l$ is the leftmost nonzero component and $M^r$ is the rightmost nonzero component, then the {\em width} of $M^\bullet$ is defined to be $r-l+1$. If $M^\bullet$ is zero, then the width of $M^\bullet$ is defined to be $0$.

\begin{proposition}\label{prop structure constans in DH(A)}
Let $A_1^\bullet=(A_1^i, 0)_{i\in\mathbb{Z}}, A_2^\bullet=(A_2^i, 0)_{i\in\mathbb{Z}}\in\cc_0(\ca)$. Then in $\cd\ch_0(\ca)$ we have that
$$[A_1^\bullet][A_2^\bullet]=\sum\limits_{X^{\bullet}=(X^i,0)_{i\in\mathbb{Z}}\in\Iso(\cd_0(\ca))}\sum_{\begin{array}{c}I^i,M^i,N^i \in\Iso(\ca)\\i\in\mathbb{Z}\end{array}}$$
$$\left(\prod_{i\in\mathbb{Z}}g^{A_1^i}_{I^iM^i}g^{X^i}_{M^iN^i}g^{A_2^i}_{N^iI^{i-1}}
\frac{a_{N^i}a_{M^i}a_{I^i}}{a_{A_1^i}a_{A_2^i}}\frac{\prod_{k\geq 2}\langle\widehat{A_2^{i+k}},\widehat{A_1^{i}}\rangle^{(-1)^k}}{\langle\widehat{N^{i+1}},\widehat{M^i}\rangle}\right)[X^\bullet].
$$
\end{proposition}

\begin{proof}
Since $A_1^\bullet$ and $A_2^\bullet$ are both bounded complexes with zero differentials, there exist $l, r\in\mathbb{Z}$ such that
$A_1^i=0$ and $A_2^i=0$ for any $i<l$ or $i>r$.  And it suffices to show the proposition in the case of $l=0$ since we can use the shift functor. So we can prove it by induction on the width of the complexes. By Lemma \ref{lemma multiplication of large item product small item} one can obtain that
$[A_i^\bullet]=Z_{A_{i}^{r}}^{[r]}Z_{A_{i}^{r-1}}^{[r-1]}\cdots Z_{A_{i}^{0}}^{[0]}$ for $i=1, 2$.

This is clear if $r=0$. Assume that $r=1$, then we have
\begin{eqnarray*}
&&[A_1^\bullet][A_2^\bullet]\\
&=&Z_{A_{1}^{1}}^{[1]}Z_{A_{1}^{0}}^{[0]}Z_{A_{2}^{1}}^{[1]}Z_{A_{2}^{0}}^{[0]}\\
&=&\sum_{M^0, N^1\in\Iso(\ca)}\gamma_{A_2^1A_1^0}^{M^0N^1}\frac{1}{\langle\widehat{N^1}, \widehat{M^0}\rangle}Z_{A_{1}^{1}}^{[1]}Z_{N^{1}}^{[1]}Z_{M^0}^{[0]}Z_{A_{2}^{0}}^{[0]}\\
&=&\sum_{X^0,X^1\in\Iso(\ca)}\left(\sum_{M^0, N^1, I^0\in\Iso(\ca)}g_{A_1^1N^1}^{X^1}g_{N^1I^0}^{A_2^1}g_{I^0M^0}^{A_1^0}g_{M^0A_2^0}^{X^0}
\frac{a_{I^0}a_{M^0}a_{N^1}}{a_{A_1^0}a_{A_2^1}}\frac{1}{\langle\widehat{N^1}, \widehat{M^0}\rangle}\right)Z_{X^{1}}^{[1]}Z_{X^{0}}^{[0]}\\
&=&\sum_{X^0,X^1\in\Iso(\ca)}\sum_{\begin{array}{c}N^0\cong A_2^0,M^0,I^0\\
N^1 ,M^1\cong A_1^1,\\
I^1\cong0\in\Iso(\ca)\end{array}}\left(\prod_{i=0,1}g_{I^iM^i}^{A_1^i}g_{M^iA_2^i}^{X^i}g_{N^iI^{i-1}}^{A_2^i}
\frac{a_{N^i}a_{M^i}a_{I^i}}{a_{A_1^i}a_{A_2^i}}\frac{1}{\langle\widehat{N^1}, \widehat{M^0}\rangle}\right)Z_{X^{1}}^{[1]}Z_{X^{0}}^{[0]}.
\end{eqnarray*}
Since $\prod_{i\in\mathbb{Z}}\prod_{k\geq 2}\langle \widehat{A_2^{i+k}},\widehat{A_1^i}\rangle=1$, it is not hard to get the equality as required.

If $r>1$, then by Proposition \ref{proposition derived hall algebra} one can obtain that
\begin{eqnarray*}
&&[A_1^\bullet][A_2^\bullet]\\
&=&Z_{A_{1}^{r}}^{[r]}Z_{A_{1}^{r-1}}^{[r-1]}\cdots Z_{A_{1}^{0}}^{[0]}Z_{A_{2}^{r}}^{[r]}Z_{A_{2}^{r-1}}^{[r-1]}\cdots Z_{A_{2}^{0}}^{[0]}\\
&=&\langle \widehat{A_2^r},\widehat{A_1^0}\rangle^{(-1)^r}\langle \widehat{A_2^r},\widehat{A_1^1}\rangle^{(-1)^{r-1}}\cdots\langle \widehat{A_2^r},\widehat{A_1^{r-2}}\rangle^{(-1)^2} Z_{A_1^r}^{[r]}Z_{A_1^{r-1}}^{[r-1]}Z_{A_2^r}^{[r]}Z_{A_1^{r-2}}^{[r-2]}\cdots Z_{A_1^0}^{[0]}Z_{A_2^{r-1}}^{[r-1]}\cdots Z_{A_2^{0}}^{[0]}\\
&=&\langle \widehat{A_2^r},\widehat{A_1^0}\rangle^{(-1)^r}\langle \widehat{A_2^r},\widehat{A_1^1}\rangle^{(-1)^{r-1}}\cdots\langle \widehat{A_2^r},\widehat{A_1^{r-2}}\rangle^{(-1)^2} \\
&&\sum_{M^{r-1},N^r\in\Iso(\ca)}\gamma_{A_2^rA_1^{r-1}}^{M^{r-1}N^r}\frac{1}{\langle \widehat{N^r},\widehat{M^{r-1}}\rangle}Z_{A_1^r}^{[r]}Z_{N^r}^{[r]}Z_{M^{r-1}}^{[r-1]}Z_{A_1^{r-2}}^{[r-2]}\cdots Z_{A_1^0}^{[0]}Z_{A_2^{r-1}}^{[r-1]}\cdots Z_{A_2^{0}}^{[0]}\\
&=&\langle \widehat{A_2^r},\widehat{A_1^0}\rangle^{(-1)^r}\langle \widehat{A_2^r},\widehat{A_1^1}\rangle^{(-1)^{r-1}}\cdots\langle \widehat{A_2^r},\widehat{A_1^{r-2}}\rangle^{(-1)^2} \sum_{N^r, M^{r-1},I^{r-1}\in\Iso(\ca)}\\
&&\left(g_{N^rI^{r-1}}^{A_2^r}g_{I^{r-1}M^{r-1}}^{A_1^{r-1}}
\frac{a_{N^{r}}a_{M^{r-1}}a_{I^{r-1}}}{a_{A_1^{r-1}}a_{A_2^r}}\frac{1}{\langle \widehat{N^r},\widehat{M^{r-1}}\rangle}\right)Z_{A_1^r}^{[r]}Z_{N^r}^{[r]}Z_{M^{r-1}}^{[r-1]}Z_{A_1^{r-2}}^{[r-2]}\cdots Z_{A_1^0}^{[0]}Z_{A_2^{r-1}}^{[r-1]}\cdots Z_{A_2^{0}}^{[0]}\\
&=&\langle \widehat{A_2^r},\widehat{A_1^0}\rangle^{(-1)^r}\langle \widehat{A_2^r},\widehat{A_1^1}\rangle^{(-1)^{r-1}}\cdots\langle \widehat{A_2^r},\widehat{A_1^{r-2}}\rangle^{(-1)^2}\sum_{X^r\in\Iso(\ca)}\left(\sum_{N^r, M^{r-1},I^{r-1}\in\Iso(\ca)}\right.\\
&&\left.g_{A_1^rN^r}^{X^r}g_{N^rI^{r-1}}^{A_2^r}g_{I^{r-1}M^{r-1}}^{A_1^{r-1}}\frac{a_{N^{r}}a_{M^{r-1}}a_{I^{r-1}}}{a_{A_1^{r-1}}a_{A_2^r}}\frac{1}{\langle \widehat{N^r},\widehat{M^{r-1}}\rangle}\right)Z_{X^r}^{[r]}Z_{M^{r-1}}^{[r-1]}Z_{A_1^{r-2}}^{[r-2]}\cdots Z_{A_1^0}^{[0]}Z_{A_2^{r-1}}^{[r-1]}\cdots Z_{A_2^{0}}^{[0]}.
\end{eqnarray*}
Since $M^r\cong A_1^r, I^r\cong0$, one can easily obtain that $a_{M^r}=a_{A_1^r}$ and $g_{I^rM^r}^{A_1^r}=1$.
By induction on the width of the complexes, we have that
$$Z_{M^{r-1}}^{[r-1]}Z_{A_1^{r-2}}^{[r-2]}\cdots Z_{A_1^0}^{[0]}Z_{A_2^{r-1}}^{[r-1]}\cdots Z_{A_2^{0}}^{[0]}=\sum\limits_{X^{\bullet}=(X^i,0)_{i\in\mathbb{Z}}\in\Iso(\cd_0(\ca))}$$
$$\left(\sum_{{\small\begin{array}{c}I^i,M^i,N^i\in\Iso(\mathcal{A})\\ i\in\mathbb{Z}\end{array}}}\left(\prod_{i\in\mathbb{Z}}
g^{A_1^i}_{I^iM^i}g^{X^i}_{M^iN^i}g^{A_2^i}_{N^iI^{i-1}}
\frac{a_{N^i}a_{M^i}a_{I^i}}{a_{A_1^i}a_{A_2^i}}\frac{\prod_{k\geq 2}\langle\widehat{A_2^{i+k}},\widehat{A_1^{i}}\rangle^{(-1)^k}}
{\langle\widehat{N^{i+1}},\widehat{M^i}\rangle}\right)\right)[X^\bullet],
$$\\
where $X^i=0$ if $i<0$ or $i>r-1$. In fact, we have that \begin{eqnarray*}
&&\sum_{\begin{array}{c}I^i,M^i,N^i\in\Iso(\mathcal{A})\\ i\in\mathbb{Z}\end{array}}\left(\prod_{i\in\mathbb{Z}}g^{A_1^i}_{I^iM^i}g^{X^i}_{M^iN^i}g^{A_2^i}_{N^iI^{i-1}}
\frac{a_{N^i}a_{M^i}a_{I^i}}{a_{A_1^i}a_{A_2^i}}\frac{\prod_{k\geq 2}\langle\widehat{A_2^{i+k}},\widehat{A_1^{i}}\rangle^{(-1)^k}}
{\langle\widehat{N^{i+1}},\widehat{M^i}\rangle}\right)\\
&=&\sum_{\begin{array}{c}I^i,M^i,N^i\in\Iso(\mathcal{A})\\ i\in\mathbb{Z}\end{array}}g^{X^{r-1}}_{M^{r-1}N^{r-1}}g^{A_2^{r-1}}_{N^{r-1}I^{r-2}}
%g^{A_1^{r-2}}_{I^{r-2M^{r-2}}}
\cdots g^{A_1^0}_{I^0M^0}g^{X^0}_{M^0A_2^0}
\frac{a_{N^{r-1}}a_{I^{r-2}}a_{M^{r-2}}\cdots a_{M^0}}{a_{A_1^{r-2}}a_{A_2^{r-1}}\cdots a_{A_2^1}a_{A_1^0}}\\
&&\frac{\langle \widehat{A_2^{r-1}},\widehat{A_1^{0}}\rangle^{(-1)^{r-1}}
\cdots\langle\widehat{A_2^{r-1}},\widehat{A_1^{r-3}}\rangle^{(-1)^{2}}
\cdots\langle\widehat{A_2^{2}},\widehat{A_1^{0}}\rangle^{(-1)^{2}}}{\langle \widehat{N^{r-1}},\widehat{M^{r-2}}\rangle\langle \widehat{N^{r-2}},\widehat{M^{r-3}}\rangle\cdots
\langle\widehat{N^{1}},\widehat{M^{0}}\rangle}.
\end{eqnarray*}
Thus one can get the multiplication as required by the isomorphism $N^0\cong A_2^0$.
\end{proof}

\subsection{Derived Hall numbers for odd periodic complexes}
Similar to the work of \cite{W18}, we can get the derived Hall number for odd periodic complexes.
\begin{lemma}[\cite{LinP2021}, Proposition5.5]\label{proposition Hall number of MHA and DHA}
Let $A_1^\bullet$ and $A_2^\bullet$ be $t$-periodic complexes in $\cc_t(\ca)$ with zero differentials. Then we have the following isomorphism
$$\Ext_{\cc_t(\ca)}^1(A_1^\bullet, A_2^\bullet)\cong\Hom_{\cd_t(\ca)}(A_1^\bullet, A_2^\bullet[1]).$$
\end{lemma}

\begin{proposition}\label{prop multiplication of complex in DHt(A)}
Let $A_1^\bullet=(A_1^i, 0)_{i\in\mathbb{Z}_t}, A_2^\bullet=(A_2^i, 0)_{i\in\mathbb{Z}_t}\in\cc_t(\ca)$. Then in $\cd_t(\ca)$ we have that
\begin{eqnarray*}
\prod\limits_{i=0}^{t-1}|\Hom_{\cd_t(\ca)}(A_1^\bullet[i],A_2^\bullet)|^{(-1)^i}&=&
\prod\limits_{i=0}^{t-1}|\Hom_{\ca}(A_1^i,A_2^i)||\Ext_{\ca}^1(A_1^i,A_2^i)|\left(\prod\limits_{k=1}^{t-1}\langle A_1^{i+k},A_2^i\rangle^{(-1)^{k}}\right)
\end{eqnarray*}
\end{proposition}
\begin{proof}
If $t=1$, it is apparently implied by the isomorphism $$\Hom_{\cd_1(\ca)}(A_1^\bullet,A_2^\bullet)\cong\Hom_{\ca}(A_1,A_2)\oplus\Ext^1_{\ca}(A_1,A_2).$$

If $t\geq 3$, we have
$$\Hom_{\cd_t(\ca)}(A_1^\bullet,A_2^\bullet)\cong\bigoplus\limits_{j=0}^{t-1}\left(\Hom_{\ca}(A_1^j,A_2^j)\bigoplus\Ext^1_{\ca}(A_1^j,A_2^{j-1})\right).$$
Therefore, one can obtain that
\begin{eqnarray*}
&\prod\limits_{i=0}^{t-1}|\Hom_{\cd_t(\ca)}(A_1^\bullet[i],A_2^\bullet)|^{(-1)^i}\hspace{10cm}\\
=&|\Hom_{\ca}(A_1^0,A_2^0)||\Hom_{\ca}(A_1^1,A_2^1)|\cdots|\Hom_{\ca}(A_1^{t-1},A_2^{t-1})||\Ext^1_{\ca}(A_1^1,A_2^0)||\Ext^1_{\ca}(A_1^2,A_2^1)|\cdots\\
&|\Ext^1_{\ca}(A_1^0,A_2^{t-1})||\Hom_{\ca}(A_1^1,A_2^0)|^{-1}|\Hom_{\ca}(A_1^2,A_2^1)|^{-1}\cdots|\Hom_{\ca}(A_1^{0},A_2^{t-1})|^{-1}|\Ext^1_{\ca}(A_1^2,A_2^0)|^{-1}\\
&|\Ext^1_{\ca}(A_1^3,A_2^1)|^{-1}\cdots|\Ext^1_{\ca}(A_1^1,A_2^{t-1})|^{-1}|\Hom_{\ca}(A_1^2,A_2^0)||\Hom_{\ca}(A_1^3,A_2^1)|\cdots|\Hom_{\ca}(A_1^{1},A_2^{t-1})|\\
&|\Ext^1_{\ca}(A_1^3,A_2^0)||\Ext^1_{\ca}(A_1^4,A_2^1)|\cdots|\Ext^1_{\ca}(A_1^2,A_2^{t-1})|\cdots\hspace{7cm}\\
&|\Hom_{\ca}(A_1^{t-2},A_2^0)|^{-1}|\Hom_{\ca}(A_1^{t-1},A_2^1)|^{-1}|\Hom_{\ca}(A_1^{0},A_2^2)|^{-1}\cdots|\Hom_{\ca}(A_1^{t-3},A_2^{t-1})|^{-1}\\
&|\Ext^1_{\ca}(A_1^{t-1},A_2^0)|^{-1}|\Ext^1_{\ca}(A_1^0,A_2^1)|^{-1}|\Ext^1_{\ca}(A_1^1,A_2^2)|^{-1}\cdots|\Ext^1_{\ca}(A_1^{t-2},A_2^{t-1})|^{-1}\\
&|\Hom_{\ca}(A_1^{t-1},A_2^0)||\Hom_{\ca}(A_1^{0},A_2^1)|\cdots|\Hom_{\ca}(A_1^{t-2},A_2^{t-1})||\Ext^1_{\ca}(A_1^{0},A_2^0)|\cdots|\Ext^1_{\ca}(A_1^{t-1},A_2^{t-1})|\\
=&|\Hom_{\ca}(A_1^0,A_2^0)|\cdots|\Hom_{\ca}(A_1^{t-1},A_2^{t-1})||\Ext^1_{\ca}(A_1^{0},A_2^0)|\cdots|\Ext^1_{\ca}(A_1^{t-1},A_2^{t-1})|\hspace{2cm}\\
&\frac{1}{\langle A_1^1,A_2^0\rangle\langle A_1^2,A_2^1\rangle\cdots\langle A_1^0,A_2^{t-1}\rangle}\langle A_1^2,A_2^0\rangle\langle A_1^3,A_2^1\rangle\cdots\langle A_1^1,A_2^{t-1}\rangle\cdots\langle A_1^{t-1},A_2^0\rangle\langle A_1^0,A_2^1\rangle\cdots\langle A_1^{t-2},A_2^{t-1}\rangle\\
=&\prod\limits_{i=0}^{t-1}|\Hom_{\ca}(A_1^i,A_2^i)||\Ext_{\ca}^1(A_1^i,A_2^i)|\left(\prod\limits_{k=1}^{t-1}\langle A_1^{i+k},A_2^i\rangle^{(-1)^{k}}\right).\hspace{6cm}
\end{eqnarray*}
\end{proof}
In the following, for each $i\in\mathbb{Z}_t$, let $A_1^i, A_2^i,X^i,S^i, M^i, N^i$ be fixed objects in $\mathcal{A}$. To simplify the notation, we set $$(f_1^\bullet, g_1^\bullet, f_2^\bullet, g_2^\bullet,f_3^\bullet, g_3^\bullet)=(f_1^0, g_1^0, f_2^0, g_2^0,f_3^0, g_3^0,\cdots, f_1^{t-1}, g_1^{t-1}, f_2^{t-1}, g_2^{t-1},f_3^{t-1}, g_3^{t-1}),$$
and $T=\{(f_1^\bullet, g_1^\bullet, f_2^\bullet, g_2^\bullet,f_3^\bullet, g_3^\bullet)|0\rightarrow S^i\xrightarrow{f_1^i} A_2^i\xrightarrow{g_1^i} N^i\rightarrow 0, 0\rightarrow N^i\xrightarrow{f_2^i} X^i\xrightarrow{g_2^i} M^i\rightarrow 0~\text{and}~0\rightarrow M^i\xrightarrow{f_3^i} A_1^i\xrightarrow{g_3^i} S^{i+1}\rightarrow 0~\text{are all short exact sequences in}~\mathcal{A}~\text{for any}~i\in\mathbb{Z}_t\}$. For a $3t$-periodic complex
\[\xymatrix{
A_2^0\ar[r]^{m^0}&X^0\ar[r]^{n^0}&A_1^0\ar[r]^{h^0}&A_2^1\ar[r]^{m^1}&X^1\ar[r]^{n^1}&A_1^1\ar[d]^{h^1}\\
A_1^{t-1}\ar[u]^{h^{t-1}}&&&&&A_2^2\ar[d]^{m^2}\\
X^{t-1}\ar[u]^{n^{t-1}}&&&&&X^2\ar[d]^{n^2}\\
A_2^{t-1}\ar[u]^{m^{t-1}}&\cdots\ar[l]_{h^{t-2}}&\cdots&\cdots&\cdots&A_1^2\ar[l]_{h^2},
}\]
set $(m^\bullet, n^\bullet, h^\bullet):=(m^0, n^0, h^0, m^1, n^1, h^1,\cdots,m^{t-1}, n^{t-1}, h^{t-1})$.

\begin{proposition}\label{prop epic}
For each $i\in\mathbb{Z}_t$, let $A_1^i, A_2^i, X^i, S^i, M^i, N^i$ be fixed objects in $\mathcal{A}$. Let $T=\{(f_1^\bullet, g_1^\bullet, f_2^\bullet, g_2^\bullet,f_3^\bullet, g_3^\bullet)\}$ and $(m^\bullet, n^\bullet, h^\bullet)$ be the set and the $3t$-periodic complex respectively as above. Set $G_{M^\bullet, N^\bullet}:=\{(m^\bullet, n^\bullet, h^\bullet)|\Im(h^i)\cong S^{i+1}, \Ker(h^i)\cong M^i,\Coker(h^i)\cong N^{i+1}\}$. Then there exists a surjective map
$$
\Phi: T\twoheadrightarrow G_{M^\bullet, N^\bullet}
$$
and then we have that $$|T|=\prod_{i\in\mathbb{Z}_t}g_{S^{i+1}M^i}^{A_1^i}g_{M^i N^i}^{X^i}g_{N^iS^i}^{A_2^i}a^2_{M^i}a^2_{N^i}a^2_{S^i}$$
and
$$|G_{M^\bullet, N^\bullet}|=\prod_{i\in\mathbb{Z}_t}g_{S^{i+1}M^i}^{A_1^i}g_{M^i N^i}^{X^i}g_{N^iS^i}^{A_2^i}a_{M^i}a_{N^i}a_{S^i}.$$
\end{proposition}

\begin{proof}
By the proof of Lemma 3.1 in \cite{P97}, we know that
the cardinality of the set $\{(f, g)\in\Hom_{\mathcal{A}}(B, C)\times\Hom_{\mathcal{A}}(C, A)|0\rightarrow B\xrightarrow{f}C\xrightarrow{g}A\rightarrow 0~\text{exact}\}$ is $\frac{|\Ext^1_{\mathcal{A}}(A,B)_C|a_C}{|\Hom_{\mathcal{A}}(A,B)|}$. Thus by the homological formula
\[g^C_{AB}=\frac{|\Ext_{\ca}^1(A, B)_C|}{|\Hom_{\ca}(A, B)|}\frac{a_C}{a_Aa_B}\]
it follows that
\[|T|=\prod_{i\in\mathbb{Z}_t}g_{S^{i+1}M^i}^{A_1^i}g_{M^i N^i}^{X^i}g_{N^iS^i}^{A_2^i}a^2_{M^i}a^2_{N^i}a^2_{S^i}\] as required.

Next one can easily get the map $\Phi$ by setting
\[(f_1^\bullet, g_1^\bullet, f_2^\bullet, g_2^\bullet,f_3^\bullet, g_3^\bullet)\mapsto (f_2^0g_1^0, f_3^0g_2^0, f_1^1g_3^0,\cdots,f_2^{t-1}g_1^{t-1}, f_3^{t-1}g_2^{t-1}, f_1^{0}g_3^{t-1})
\] which is obviously surjective, and define the free action of $\aut(M^\bullet)\times\aut(N^\bullet)\times\aut(S^\bullet):=
\aut(M^0)\times\aut(N^0)\times\aut(S^0)\times\cdots\times\aut(M^{t-1})\times\aut(N^{t-1})\times\aut(S^{t-1})$
on the set $T$ by setting
\[(l_{M^\bullet}, l_{N^\bullet}, l_{S^\bullet})(f_1^\bullet, g_1^\bullet, f_2^\bullet, g_2^\bullet,f_3^\bullet, g_3^\bullet):=(f_1^\bullet l^{-1}_{S^\bullet}, l_{N^\bullet}g_1^\bullet, f_2^\bullet l^{-1}_{N^\bullet}, l_{M^\bullet}g_2^\bullet, f_3^\bullet l^{-1}_{M^\bullet}, l_{S^\bullet}g_3^\bullet)
\]
which denotes the morphisms
$$\left(f_1^0l^{-1}_{S^0},l_{N^0}g_1^0,f_2^0l_{N^0}^{-1},l_{M^0}g_2^0,f_3^0l_{M^0}^{-1},l_{S^1}g_3^0,
\cdots,f_{1}^{t-1}l_{S^{t-1}}^{-1},l_{N^{t-1}}g_1^{t-1},f_2^{t-1}l_{N^{t-1}}^{-1},l_{M^{t-1}}g_2^{t-1},
f_{3}^{t-1}l_{M^{t-1}}^{-1},l_{S^0}g_3^{t-1}\right).$$
It follows that
\[|G_{M^\bullet, N^\bullet}|=\frac{|T|}{\prod_{i\in\mathbb{Z}_t}a_{M^i}a_{N^i}a_{S^i}},\]which implies the formula as required.
\end{proof}

In the following, we simply write $\aut(X^\bullet)$ to represent $\aut(X^0)\times\aut(X^1)\times\cdots\times\aut(X^{t-1})$. The action of $\aut(X^\bullet)$ on $G_{M^\bullet, N^\bullet}$ induces the orbit space
\[
\mathcal{G}_{M^\bullet,N^\bullet}=\{(m^\bullet, n^\bullet, h^\bullet)^{\wedge}|(m^\bullet, n^\bullet, h^\bullet)\in G_{M^\bullet,N^\bullet}\},\]
where
\[(m^\bullet, n^\bullet, h^\bullet)^{\wedge}:=\{(l_{X^0}m^0, n^0l^{-1}_{X^0}, h^0,\cdots, l_{X^{t-1}}m^{t-1}, n^{t-1}l^{-1}_{X^{t-1}}, h^{t-1})|(l_{X^0}, l_{X^1},\cdots, l_{X^{t-1}})\in\aut(X^\bullet)\}.
\]
Apparently, for any $(m^\bullet, n^\bullet, h^\bullet)^{\wedge}\in\mathcal{G}_{M^\bullet,N^\bullet}$, the stable subgroup is isomorphic to
$$\Hom_{\mathcal{A}}(M^0,N^0)\times\Hom_{\mathcal{A}}(M^1,N^1)\times\cdots\times\Hom_{\mathcal{A}}(M^{t-1},N^{t-1}),$$
so $$|\mathcal{G}_{M^\bullet,N^\bullet}|=
\frac{|G_{M^\bullet,N^\bullet}||\Hom_{\mathcal{A}}(M^0,N^0)\times\cdots\times\Hom_{\mathcal{A}}(M^{t-1},N^{t-1})|}{|\aut(X^\bullet)|}.$$

\begin{proposition}\label{thm onto mapping}
Let $A_1^\bullet=(A_1^i, 0)_{i\in\mathbb{Z}_t}, A_2^\bullet=(A_2^i, 0)_{i\in\mathbb{Z}_t}, X^\bullet=(X^i, 0)_{i\in\mathbb{Z}_t}\in\cc_t(\ca)$ and $S^i\in\mathcal{A}$ for any $i\in\mathbb{Z}_t$. Write $$\varepsilon:=\bigsqcup\limits_{\begin{array}{c}C^\bullet\in\Iso(C_t(\ca))\\\mathrm{H}(C^\bullet)\cong X^\bullet\in\Iso(D_t(\ca))\\\Im(d^i_{C^\bullet})\cong S^{i+1}, i\in\mathbb{Z}_t\end{array}}\Ext^1_{C_{t}(\ca)}(A_1^\bullet,A_2^\bullet)_{C^\bullet}$$
and
$$\mathcal{G}:=\bigsqcup\limits_{\begin{array}{c}M^i\in\Iso(\ca),N^i\in\Iso(\ca)\\i\in\mathbb{Z}_t\end{array}}\mathcal{G}_{M^\bullet,N^\bullet}.$$ Then there is a surjective map $\Phi:\varepsilon\twoheadrightarrow\mathcal{G}$.
\end{proposition}
\begin{proof}
For any $[\xi]\in\varepsilon$, assume that it is a short exact sequence of the following form
\[\xymatrix{
A_2^0\ar@{>->}[d]_{k^0}\ar[r]^{0}&A_2^1\ar[r]^{0}\ar@{>->}[d]_{k^1}&\cdots\ar[r]^{0}&A_2^{t-1}\ar@{>->}[d]_{k^{t-1}}\ar[r]^0&A_2^0\ar@{>->}[d]_{k^0}\\
C^0\ar@{->>}[d]_{p^0}\ar[r]^{d^0_{C^\bullet}}&C^1\ar[r]^{d^1_{C^\bullet}}\ar@{->>}[d]_{p^1}&\cdots\ar[r]^{d^{t-2}_{C^\bullet}}&C^{t-1}\ar@{->>}[d]_{p^{t-1}}\ar[r]^{d^{t-1}_{C^\bullet}}&C^0\ar@{->>}[d]_{p^0}\\
A_1^0\ar[r]^{0}&A_1^1\ar[r]^{0}&\cdots\ar[r]^{0}&A_1^{t-1}\ar[r]^{0}&A_1^0.
}
\]
And then, for each $i\in\mathbb{Z}_t$, there exists unique $h^i: A_1^i\rightarrow A_2^{i+1}$ such that $d^i_{C^\bullet}=k^{i+1}h^ip^i$ with $\Im(h^i)\cong \Im(d^i_{C^\bullet})\cong S^{i+1}$.

For each $h^{i-1}: A_1^{i-1}\rightarrow A_2^{i}$, denote by $0\rightarrow\Im(h^{i-1})\xrightarrow{f_1^i}A_2^i\xrightarrow{g_1^{i}}\Coker(h^{i-1})\rightarrow 0$ the canonical short exact sequence, and then form the push-out of $(k^i, g_1^i)$ and obtain the following commutative diagram
\[
\xymatrix{
\Im(h^{i-1})\ar@{>-->}[d]_{f_1^{i}}\ar@{=}[r]&\Im(h^{i-1})\ar[d]_{k^if_1^{i}}\\
A_2^i\ar@{-->>}[d]_{g_1^i}\ar[r]^{k^i}&C^i\ar[r]^{p^i}\ar@{-->}_{s^i}[d]&A_1^i\ar@{=}[d]\\
\Coker(h^{i-1})\ar@{-->}[r]^{u^i}&Z^i\ar@{-->}[r]^{v^i}&A_1^i.
}
\]
Similarly, we form the pull-back of $(v^i, f_3^i)$, where $f_3^i: \Ker(h^i)\rightarrow A_1^i$. So we get the commutative diagram
\[
\xymatrix{
\Im(h^{i-1})\ar[d]_{f_1^{i}}\ar@{=}[r]&\Im(h^{i-1})\ar[d]_{k^if_1^{i}}\\
A_2^i\ar[d]_{g_1^i}\ar[r]^{k^i}&C^i\ar[r]^{p^i}\ar@{-->}_{s^i}[d]&A_1^i\ar@{=}[d]\\
\Coker(h^{i-1})\ar@{=}[d]\ar@{-->}[r]^{u^i}&Z^i\ar@{-->}[r]^{v^i}&A_1^i\\
\Coker(h^{i-1})\ar@{-->}[r]^{f_2^i}&X^i\ar@{-->}[r]^{g_2^i}\ar@{-->}[u]^{t^i}&\Ker(h^i)\ar[u]^{f_3^i},
}
\]
in which the compositions $f_2^ig_1^i$ and $f_3^ig_2^i$ are both unique, and denote them by $m^i$ and $n^i$ independently. Thus we get a $3t$-periodic complex
\[\xymatrix{
A_2^0\ar[r]^{m^0}&X^0\ar[r]^{n^0}&A_1^0\ar[r]^{h^0}&A_2^1\ar[r]^{m^1}&X^1\ar[r]^{n^1}&A_1^1\ar[d]^{h^1}\\
A_1^{t-1}\ar[u]^{h^{t-1}}&&&&&A_2^2\ar[d]^{m^2}\\
X^{t-1}\ar[u]^{n^{t-1}}&&&&&X^2\ar[d]^{n^2}\\
A_2^{t-1}\ar[u]^{m^{t-1}}&\cdots\ar[l]_{h^t-2}&\cdots&\cdots&\cdots&A_1^2\ar[l]_{h^2},
}\]
with $\Im(h^i)\cong S^{i+1}$. So the map $\Phi$ is well-defined.

On the other hand, for any $\omega^{\wedge}\in\mathcal{G}_{M^\bullet,N^\bullet}$, assume that it belongs to the orbit of the above $3t$-periodic acyclic complex under the action of $\aut(X^\bullet)$.
Since $\mathcal{A}$ is hereditary, induced by the standard short exact sequence $0\rightarrow \Coker(h^{i-1})\xrightarrow{f_2^i}X^i\xrightarrow{g_2^i}\Ker(h^i)\rightarrow0$
and the monomorphism $f_3^i: \Ker(h^i)\rightarrow A_1^i$ we have the following commutative diagram
\[
\xymatrix{
\Coker(h^{i-1})\ar@{>->}[r]^{f_2^i}\ar@{=}[d]&X^i\ar[r]^{g_2^i}\ar@{-->}[d]_{t^i}&\Ker(h^i)\ar@{>->}[d]_{f_3^i}\\
\Coker(h^{i-1})\ar@{-->}[r]^{u^i}&Z^i\ar@{-->}[r]^{v^i}&A_1^i\\
A_2^i\ar@{->>}[u]^{g_1^i}\ar@{>-->}[r]^{k^i}&C^i\ar@{-->>}[r]^{p^i}\ar@{-->>}[u]^{s^i}&A_1^i\ar@{=}[u]
}\]
in $\mathcal{A}$ such that all the rows are exact. By setting $d_{C^\bullet}^i:=k^{i+1}h^ip^i$, we get a $t$-periodic complex $(C^\bullet, d_{C^\bullet})$ with $\Im(d^i_{C^\bullet})\cong\Im(h^i)\cong S^{i+1}$ and a short exact sequence $\eta$ of $t$-periodic complexes as follows
$$0\rightarrow A_2^\bullet\xrightarrow{k^\bullet}C^\bullet\xrightarrow{p^\bullet}A_1^\bullet\rightarrow0$$
which satisfies $\mathrm{H}(C^\bullet)\cong X^\bullet$ and $\Phi([\eta])=\omega^{\wedge}$.
\end{proof}

\begin{proposition}\label{Prop cardinality of hom}
For some fixed $S^i\in\Iso(\mathcal{A})(i\in\mathbb{Z}_t)$, the cardinality of the set $\varepsilon$ defined in Proposition \ref{thm onto mapping} is $$\sum\limits_{\begin{array}{c}M^j,N^j\in\Iso(\ca)\\j\in\mathbb{Z}_t\end{array}}
\left(\prod\limits_{i\in\mathbb{Z}_t}\frac{g_{S^{i+1}M^i}^{A_1^i}g_{M^iN^i}^{X^i}g_{N^iS^i}^{A_2^i}a_{M^i}a_{N^i}a_{S^i}|\Hom_\ca(A_1^i,A_2^i)|}
{a_{X^i}\langle A_1^i,S^i\rangle\langle S^{i+1},N^i\rangle}\right).$$
\end{proposition}

\begin{proof}
For any $\omega^{\wedge}\in\mathcal{G}_{M^\bullet,N^\bullet}$ and each $i\in\mathbb{Z}_t$,
by the proof of Proposition \ref{thm onto mapping}, we get the following two exact sequences\\
$0\rightarrow\Hom_{\mathcal{A}}(\Im(h^i), \Coker(h^{i-1}))\rightarrow\Hom_{\mathcal{A}}(A_1^i, \Coker(h^{i-1}))\rightarrow\Hom_{\mathcal{A}}(\Ker(h^i), \Coker(h^{i-1}))\rightarrow \Ext^1_{\mathcal{A}}(\Im(h^i), \Coker(h^{i-1}))\rightarrow\Ext^1_{\mathcal{A}}(A_1^i, \Coker(h^{i-1}))\xrightarrow{(f_3^i)^*}\Ext^1_{\mathcal{A}}(\Ker(h^i), \Coker(h^{i-1}))\rightarrow0$\\
and\\
$0\rightarrow\Hom_{\mathcal{A}}(A_1^i, \Im(h^{i-1}))\rightarrow \Hom_{\mathcal{A}}(A_1^i,A_2^i)\rightarrow\Hom_{\mathcal{A}}(A_1^i,\Coker(h^{i-1}))\rightarrow\Ext^1_{\mathcal{A}}(A_1^i, \Im(h^{i-1}))\rightarrow \Ext^1_{\mathcal{A}}(A_1^i,A_2^i)\xrightarrow{(g_1^i)_*}\Ext^1_{\mathcal{A}}(A_1^i,\Coker(h^{i-1}))\rightarrow 0$, so
\begin{eqnarray*}
|\Phi^{-1}(\omega^{\wedge})|&=&\prod\limits_{i\in\mathbb{Z}_t}\left(|\Ker((f_3^i)^*)||\Ker((g_1^i)_*)|\right)\\
&=&\prod\limits_{i\in\mathbb{Z}_t}\left(
\frac{|\Hom_{\mathcal{A}}(A_1^i, \Coker(h^{i-1}))||\Ext^1_{\mathcal{A}}(\Im(h^i), \Coker(h^{i-1}))|}{|\Hom_{\mathcal{A}}(\Im(h^i), \Coker(h^{i-1}))||\Hom_{\mathcal{A}}(\Ker(h^i), \Coker(h^{i-1}))|}\right)\\
&&\prod\limits_{i\in\mathbb{Z}_t}\left(
\frac{|\Hom_{\mathcal{A}}(A_1^i,A_2^i)||\Ext^1_{\mathcal{A}}(A_1^i, \Im(h^{i-1}))|}{|\Hom_{\mathcal{A}}(A_1^i, \Im(h^{i-1}))||\Hom_{\mathcal{A}}(A_1^i,\Coker(h^{i-1}))|}\right).
\end{eqnarray*}

Since $\Im(h^i)\cong S^{i+1}, \Ker(h^i)\cong M^i,\Coker(h^i)\cong N^{i+1}$,
\begin{eqnarray*}
|\Phi^{-1}(\omega^{\wedge})|&=&\prod\limits_{i\in\mathbb{Z}_t}
\left(\frac{|\Hom_{\mathcal{A}}(A_1^i,A_2^i)|}{|\Hom_{\mathcal{A}}(\Ker(h^i), \Coker(h^{i-1}))|\langle\Im(h^i),\Coker(h^{i-1})\rangle
\langle A_1^i, \Im(h^{i-1})\rangle}\right) \\
&=&\prod\limits_{i\in\mathbb{Z}_t}
\left(\frac{|\Hom_{\mathcal{A}}(A_1^i,A_2^i)|}
{|\Hom_{\mathcal{A}}(M^i, N^i)|} \frac{1}{\langle S^{i+1},N^i\rangle
\langle A_1^i, S^i\rangle}\right).
\end{eqnarray*}

Furthermore according to the proof of Proposition \ref{thm onto mapping}, we have that
\begin{eqnarray*}
|\varepsilon|&=&\sum\limits_{\omega^{\wedge}\in\mathcal{G}}|\Phi^{-1}(\omega^{\wedge})|\\
&=&\sum\limits_{\begin{array}{c}M^j,N^j\in\Iso(\ca)\\j\in\mathbb{Z}_t\end{array}}
\left(|\Phi^{-1}(\omega^{\wedge})||\mathcal{G}_{M^\bullet,N^\bullet}|\right)\\
&=&\sum\limits_{\begin{array}{c}M^j,N^j\in\Iso(\ca)\\j\in\mathbb{Z}_t\end{array}}
\left(\prod\limits_{i\in\mathbb{Z}_t}
\left(
\frac{|\Hom_{\mathcal{A}}(A_1^i,A_2^i)|}{\langle S^{i+1},N^i\rangle
\langle A_1^i, S^i\rangle}\frac{|G_{M^\bullet,N^\bullet}|}{|\aut(X^\bullet)|}\right)\right)\\
&=&\sum\limits_{\begin{array}{c}M^j,N^j\in\Iso(\ca)\\j\in\mathbb{Z}_t\end{array}}
\left(\prod\limits_{i\in\mathbb{Z}_t}\frac{g_{S^{i+1}M^i}^{A_1^i}g_{M^iN^i}^{X^i}g_{N^iS^i}^{A_2^i}a_{M^i}a_{N^i}a_{S^i}|\Hom_\ca(A_1^i,A_2^i)|}
{a_{X^i}\langle A_1^i,S^i\rangle\langle S^{i+1},N^i\rangle}\right).\end{eqnarray*}
This completes the proof.
\end{proof}

\begin{proposition}\label{prop structure constans in DtH(A)}
Let $A_1^\bullet=(A_1^i, 0)_{i\in\mathbb{Z}_t}, A_2^\bullet=(A_2^i, 0)_{i\in\mathbb{Z}_t}\in\cc_t(\ca)$ Then in $\cd\ch_t(\ca)$ we have that
\begin{eqnarray*}
[A_1^\bullet][A_2^\bullet]&=&\sum\limits_{X^{\bullet}=(X^i,0)_{i\in\mathbb{Z}_t}\in\Iso(\cd_t(\ca))}\sum_{\begin{array}{c}
M^j,N^j,S^j\in\Iso(\mathcal{A})\\j\in \mathbb{Z}_t\end{array}}\\
&&\left(\prod_{i\in\mathbb{Z}_t}\frac{g_{S^{i+1}M^i}^{A_1^i}g_{M^iN^i}^{X^i}g_{N^iS^i}^{A_2^i}a_{M^i}a_{N^i}a_{S^i}}
{a_{X^i}\langle A_1^i,S^i\rangle\langle S^{i+1},N^i\rangle}
\sqrt{\langle A_1^i,A_2^i\rangle\prod_{k=1}^{t-1}\langle A_1^{i+k},A_2^i\rangle^{(-1)^{k+1}}}\frac{a_{X^\bullet}'}{a_{A_1^\bullet}'
a_{A_2^\bullet}'}\right)[X^\bullet],
\end{eqnarray*}
where $a_{Y^\bullet}'=a_{Y^\bullet}\{Y^\bullet, Y^\bullet\}^{\frac{1}{2}}$ for any $Y^\bullet\in\cc_t(\ca)$.
\end{proposition}

\begin{proof}
By Proposition \ref{proposition Hall number of MHA and DHA} and Proposition \ref{proposition t-periodic iso to the homology}, we know that
\[\Hom_{\mathcal{D}_t(\ca)}(A_1^\bullet, A_2^\bullet[1])_{X^\bullet[1]}\cong \bigsqcup\limits_{\begin{array}{c}C^\bullet\in\Iso(C_t(\ca))\\\mathrm{H}(C^\bullet)\cong X^\bullet\in\Iso(D_t(\ca))\end{array}}\Ext^1_{C_{t}(\ca)}(A_1^\bullet,A_2^\bullet)_{C^\bullet}\]
for any $A_1^\bullet, A_2^\bullet, \in\mathcal{C}_t(\mathcal{A})$.

Then by definition we have that
\[
[A_1^\bullet][A_2^\bullet]=\sum\limits_{X^\bullet\in\Iso(\mathcal{D}_t(\mathcal{A}))}
\frac{|\Hom_{\mathcal{D}_t(\ca)}(A_1^\bullet, A_2^\bullet[1])_{X^\bullet[1]}|}{\sqrt{\prod\limits_{i=0}^{t-1}|\Hom_{\mathcal{D}_t(\ca)}(A_1^\bullet[i],A_2^\bullet)|^{(-1)^i}}}
\frac{a_{X^\bullet}'}{a_{A_1^\bullet}'
a_{A_2^\bullet}'}
[X^\bullet],\]
and it follows from Propositions \ref{prop multiplication of complex in DHt(A)} and \ref{Prop cardinality of hom} that
\begin{eqnarray*}
&&[A_1^\bullet][A_2^\bullet]\\
&=&\sum\limits_{X^\bullet\in\Iso(\mathcal{D}_t(\mathcal{A}))}\sum\limits_{\begin{array}{c}M^j,N^j,S^j\in\Iso(\ca)\\j\in\mathbb{Z}_t\end{array}}\\
&&\left(\frac{
\prod\limits_{i\in\mathbb{Z}_t}\frac{g_{S^{i+1}M^i}^{A_1^i}g_{M^iN^i}^{X^i}g_{N^iS^i}^{A_2^i}a_{M^i}a_{N^i}a_{S^i}|\Hom_\ca(A_1^i,A_2^i)|}
{a_{X^i}\langle A_1^i,S^i\rangle\langle S^{i+1},N^i\rangle}}{\sqrt{\prod\limits_{i=0}^{t-1}|\Hom_{\ca}(A_1^i,A_2^i)||\Ext_{\ca}^1(A_1^i,A_2^i)|\left(\prod\limits_{k=1}^{t-1}\langle A_1^{i+k},A_2^i\rangle^{(-1)^{k}}\right)}}\right)\frac{a_{X^\bullet}'}{a_{A_1^\bullet}'
a_{A_2^\bullet}'}[X^\bullet]\\
&=&\sum\limits_{X^\bullet\in\Iso(\mathcal{D}_t(\mathcal{A}))}\sum_{\begin{array}{c}
M^j,N^j,S^j\in\Iso(\mathcal{A})\\j\in \mathbb{Z}_t\end{array}}\\
&&\left(\prod_{i\in\mathbb{Z}_t}\frac{g_{S^{i+1}M^i}^{A_1^i}g_{M^iN^i}^{X^i}g_{N^iS^i}^{A_2^i}a_{M^i}a_{N^i}a_{S^i}}
{a_{X^i}\langle A_1^i,S^i\rangle\langle S^{i+1},N^i\rangle}
\sqrt{\langle A_1^i,A_2^i\rangle\prod_{k=1}^{t-1}\langle A_1^{i+k},A_2^i\rangle^{(-1)^{k+1}}}\right)\frac{a_{X^\bullet}'}{a_{A_1^\bullet}'
a_{A_2^\bullet}'}[X^\bullet]
\end{eqnarray*}
\end{proof}

\subsection{Main result}\label{subsection main result}
By the definition of $\cl_t(\ca)$, Proposition \ref{prop structure constans in DH(A)} and Proposition \ref{prop structure constans in DtH(A)}, one can easily get the following natural isomorphisms between $\cl_t(\ca)$ and derived Hall algebras.

\begin{theorem}\label{theorem green to derived Hall algebras}
The associative algebra $\cl_t(\ca)$ is isomorphic to derived Hall algebra $\mathcal{D}\mathcal{H}_t(\mathcal{A})$ of $\mathcal{A}$ for any $t=0$ or $t$ is an odd positive integer.
\end{theorem}

\begin{remark}
By Theorems \ref{Thm associativity} and \ref{theorem green to derived Hall algebras}, we say that the associativity of derived Hall algebra can be deduced from the Green's formula.
\end{remark}

\end{document}